\documentclass{article}

\usepackage{url,amsmath,marvosym}
\usepackage{amsfonts}
\usepackage{graphicx}
\usepackage{epstopdf}
\usepackage{geometry}

\usepackage{fancyhdr}
\usepackage{color,marginnote,bm,hyperref}

\def\QED{\fbox{}}
\def\dd{\mathrm{d}}

\def\RR{\mathbb{R}}
\def\NN{\mathbb{N}}
\def\P{{\cal P}}
\def\I{{\cal I}}
\def\bfb{\bm{b}}
\def\bfc{\bm{c}}
\def\bfuno{\bm{1}}
\def\bfzero{\bm{0}}
\def\bfgamma{\bm{\gamma}}
\def\eps{\varepsilon}
\def\diag{\mathrm{diag}}
\def\proof{\underline{Proof}~}
\def\no{\noindent}

\newtheorem{theorem}{Theorem}
\newtheorem{corollary}{Corollary}
\newtheorem{remark}{Remark}
\newtheorem{lemma}{Lemma}

\title{A new framework for polynomial approximation to differential equations\,\footnote{
Cite as: L.\,Brugnano, G.\,Frasca-Caccia, F.\,Iavernaro, V.\,Vespri. A new framework for polynomial approximation to differential equations. {\em Adv. Comput. Math.} {\bf 48},\,76 (2022) \url{https://doi.org/10.1007/s10444-022-09992-w} }}
\author{Luigi Brugnano\footnote{Universit\`a di Firenze, Italy  (\url{luigi.brugnano@unifi.it}, orcidID:\,{0000-0002-6290-4107})}
\and Gianluca Frasca-Caccia\footnote{Universit\`a di Salerno, Italy   (\url{gfrascacaccia@unisa.it})}
\and Felice Iavernaro\footnote{Universit\`a di Bari, Italy (\url{felice.iavernaro@uniba.it}, orcidID:\,{0000-0002-9716-7370})}
\and Vincenzo Vespri\footnote{Universit\`a di Firenze, Italy (\url{vincenzo.vespri@unifi.it})}
}
\date{}

\begin{document}

\maketitle

\begin{abstract}
In this paper we discuss a framework for the polynomial approximation
to the solution of initial value problems for differential equations. The
framework is based on an expansion of the vector field along an orthonormal basis,
and relies on perturbation results for the considered problem.  Initially devised for the approximation of ordinary differential 
equations, it is here further extended  and, moreover, generalized to cope with constant delay differential 
equations. Relevant classes of Runge-Kutta methods can be derived 
within this framework.

\medskip
\no{\bf Keywords:}
  ordinary differential equations, delay differential equations, orthogonal polynomials, local Fourier expansion, polynomial approximations, Runge-Kutta methods.

\medskip
\no{\bf MCS:}  65L05, 65L03, 65L06, 65P10.
\end{abstract}

\section{Introduction}\label{sec:intro}
In this paper, we shall deal with the definition of a framework to discuss polynomial approximations to the solution of
initial value problems for ordinary differential equations (ODEs),%\nnote{\blue eq:odeivp}
\begin{equation}\label{eq:odeivp}
\dot y(t) = f(t,y(t)), \qquad t\in[t_0,T], \qquad y(t_0) = y_0\in\RR^m, 
\end{equation}
and delay differential equations (DDEs) in the form,%\nnote{\blue eq:ddeivp}
\begin{eqnarray}\label{eq:ddeivp}
\dot y(t) &=& f(t,y(t),y(t-\tau)), \qquad t\in[t_0,T], \qquad y(t_0) = y_0,\\ \nonumber
y(t)&=&\phi(t), \qquad\qquad\qquad\qquad t\in[t_0-\tau,t_0),\nonumber
\end{eqnarray}
where $\tau>0$ is a constant delay and, usually, $y_0=\phi(t_0)$. In the sequel, we shall always assume that $f$ and $\phi$ 
are suitably regular in their respective arguments. As is well known, the two problems are related in many ways but, at the
same time, have quite different features, which reflect on their numerical solution. We refer, e.g., to the comprehensive monograph
\cite{HNW08}, concerning (\ref{eq:odeivp}), and \cite{BZ03} (see also \cite{Br04}) for (\ref{eq:ddeivp}).

In more detail, in this paper we shall fully develop a novel framework for deriving numerical methods for solving (\ref{eq:odeivp}), which
is then extended to cope with (\ref{eq:ddeivp}). %Further extensions will be the subject of future investigations.

The framework we are interested in relies on a local expansion of the vector field in (\ref{eq:odeivp}) along an orthonormal basis. 
Such basis will be, in the present case, the Legendre polynomial basis $\{P_j\}_{j\ge0}$:%\nnote{\blue eq:leg}
\begin{equation}\label{eq:leg}
P_j\in\Pi_j, \qquad \int_0^1 P_i(x)P_j(x)\dd x = \delta_{ij}, \qquad i,j=0,1,\dots,
\end{equation}
where, as is usual, $\Pi_j$ is the vector space of polynomials of degree $j$, and $\delta_{ij}$ is the Kronecker symbol.
The idea is actually not new: early use of this approach are, for example, Hulme \cite{Hu72_1,Hu72_2},  Bottasso \cite{Bo97}, and 
Betsch and Steinmann \cite{BS00}; it is also at the basis of the energy-conserving class of Runge-Kutta methods named HBVMs 
\cite{BIT2010} (see also the monograph \cite{BI2016} and the review paper \cite{BI2018}). 

The approach that we shall pursue has been initially devised in \cite{BIT2012},  where the target was problem (\ref{eq:odeivp}), and its potentialities have been disclosed by using HBVMs as spectral methods in time for efficiently solving highly oscillatory problems \cite{BMR2019} and, subsequently, Hamiltonian PDEs \cite{BIMR2019}. A corresponding error analysis is given in \cite{ABI2020}. Moreover, this allows deriving a formulation of HBVMs as continuous-stage Runge-Kutta methods \cite{ABI2019,ABI2022}.

Starting from this background, in this work we carry out a complete perturbation analysis of problems (\ref{eq:odeivp}) and (\ref{eq:ddeivp}), and set up a unique and comprehensive framework to deal with the numerical solution of both problems by exploiting the same discretization procedure. In more detail, the truncated Fourier series may be interpreted as a projection of the differential problem onto a finite dimensional vector space, leading to a new, numerically easy-to-handle, differential problem. This latter may be regarded as a perturbation of the original one, so that the perturbation analysis turns out to be crucial to understand how the solutions of the two problems are related.  At the best of our knowledge, the perturbation results for problem (\ref{eq:ddeivp}) are new, and provide a powerful general tool of analysis. That the same framework may cover problems of different nature constitute, in our opinion, a specific advancement in this field,  and reveals its potentialities to deal with other classes of problems (which will be the subject of future investigations).

With this premise, the paper is organized as follows: Section~\ref{sec:ode} concerns the result pertaining to the ODE case;
Section~\ref{sec:dde} contains the corresponding results for the DDE case; 
Section~\ref{sec:numtest} contains some numerical tests for the DDE case, involving methods which are new in this setting; at last, some concluding remarks and possible developments are reported in Section~\ref{sec:conclusions}.

\section{The ODE case}\label{sec:ode} 
Without loss of generality, we shall consider problem (\ref{eq:odeivp}) in the simpler form:%\rnote{\blue eq:ode1}
\begin{equation}\label{eq:ode1} 
\dot y(t) = f(y(t)), \qquad t\in [t_0,T], \qquad y(t_0)=y_0\in\RR^m.
\end{equation}
Having fixed the mesh%\rnote{\blue eq:h}
\begin{equation}\label{eq:h}
t_n = t_0+nh, \qquad n=0,\dots,N, \qquad h = \frac{T-t_0}N,
\end{equation}
we formally set, for $n=1,\dots,N$:%\rnote{\blue eq:hsign}
\begin{equation}\label{eq:hsign}
\hat\sigma_n(ch) := y(t_{n-1}+ch) \equiv \hat\sigma(t_{n-1}+ch), \qquad c\in[0,1],
\end{equation}
the restriction of the solution of problem (\ref{eq:ode1}) to the time interval $[t_{n-1},t_n]$ (the function $\hat\sigma(t)\equiv y(t)$ is introduced for notational purposes). Consequently, $\hat\sigma_n$ satisfies the 
differential equation%\rnote{\blue eq:hsig1}
\begin{equation}\label{eq:hsig1}
\dot{\hat\sigma}_n(ch) = \sum_{j\ge 0} P_j(c)\gamma_j(\hat\sigma_n), \qquad c\in[0,1], \qquad \hat\sigma_n(0) = y(t_{n-1}),
\end{equation}
so that,%\rnote{\blue eq:hsig}
\begin{equation}\label{eq:hsig}
\hat\sigma_n(ch) = y(t_{n-1})+h\sum_{j\ge0} \int_0^c P_j(x)\dd x\,\gamma_j(\hat\sigma_n), \qquad c\in[0,1],
\end{equation}
and, by virtue of (\ref{eq:leg}),%\rnote{\blue eq:ytn}
\begin{equation}\label{eq:ytn}
y(t_n) = y(t_{n-1}) + h\gamma_0(\hat\sigma_n) \equiv \hat\sigma(t_n),
\end{equation}
where, in general, for any suitably regular function $z:[0,h]\rightarrow\RR^m$,%\rnote{\blue eq:gammaj}
\begin{equation}\label{eq:gammaj}
\gamma_j(z) := \int_0^1 P_j(\zeta)f(z(\zeta h))\dd\zeta.
\end{equation}

We now look for a piecewise polynomial approximation $\sigma(t)$, to the solution of (\ref{eq:ode1}), such that, setting for $n=1,\dots,N$,%\rnote{\blue eq:sign}
\begin{equation}\label{eq:sign}
\sigma_n(ch) \equiv \sigma(t_{n-1}+ch), \qquad c\in[0,1],
\end{equation}                          
its restriction to the time interval $[t_{n-1},t_n]$, $\sigma_n\in\Pi_s$ and satisfies the differential equation:%\note{\blue eq:sig1}
\begin{equation}\label{eq:sig1}
\dot\sigma_n(ch) = \sum_{j=0}^{s-1} P_j(c) \gamma_j(\sigma_n), \qquad c\in [0,1],\qquad
\sigma_n(0) = y_{n-1},
\end{equation}
obtained by truncating the infinite series in (\ref{eq:hsig1}) to a finite sum, with %\note{\blue eq:yn}
\begin{equation}\label{eq:yn}
y_n := \sigma_n(h) \equiv \sigma(t_n).
\end{equation}
Consequently, $\sigma_n$ can be formally written as:%\note{\blue eq:sig}
\begin{equation}\label{eq:sig}
\sigma_n(ch) = y_{n-1} + h\sum_{j=0}^{s-1} \int_0^c P_j(x)\dd x \, \gamma_j(\sigma_n), \qquad c\in[0,1],
\end{equation}
and (compare with (\ref{eq:ytn})),%\nnote{\blue eq:yn-1}
\begin{equation}\label{eq:yn-1}
y_n = y_{n-1} + h\gamma_0(\sigma_n). %, \qquad n=1,\dots,N.
\end{equation}

\subsection{Preliminary results}\label{sec:prelode}
We here provide a few preliminary results, which will be needed to derive the main ones in the following subsections. Some of them are
taken from \cite{BIT2012} but we also report them here, for sake of completeness.

\begin{theorem}\label{thm:Gh} Let $G:[0,h]\rightarrow V$, with $V$ a vector space, admit a Taylor expansion at 0. Then, for all $j\ge0$:
$$\int_0^1 P_j(\zeta) G(\zeta h)\dd\zeta = O(h^j).$$
\end{theorem}
\proof  By virtue of (\ref{eq:leg}), one has:
\begin{eqnarray*}
\int_0^1 P_j(\zeta) G(\zeta h)\dd\zeta &=& \int_0^1 P_j(\zeta) \sum_{k\ge0} \frac{G^{(k)}(0)}{k!}(\zeta h)^k\dd\zeta =  \sum_{k\ge0} \frac{G^{(k)}(0)}{k!}h^k\int_0^1P_j(\zeta)\zeta^k\dd\zeta\\
&=& \sum_{k\ge j} \frac{G^{(k)}(0)}{k!}h^k\int_0^1P_j(\zeta)\zeta^k\dd\zeta = O(h^j).\,
\QED
\end{eqnarray*}

As a straightforward consequence, setting $G(\zeta h):=f(z(\zeta h))$, the following result is proved.

\begin{corollary}\label{cor:gamj} With reference to (\ref{eq:gammaj}), one has: \,$\gamma_j(z) = O(h^j)$.
\end{corollary}

Let us denote by
\begin{equation}\label{eq:ysol}
y(t)\equiv y(t,\xi,\eta) 
\end{equation}
the solution of the problem (compare with (\ref{eq:ode1})):%\nnote{\blue eq:odexi}
\begin{equation}\label{eq:odexi}
\dot y(t) = f(y(t)), \qquad t\in[\xi,T], \qquad y(\xi)=\eta\in\RR^m.
\end{equation}
Hereafter, for sake of brevity, we may use either one of the two notations in (\ref{eq:ysol}), depending on the needs.  The following theorem contains standard perturbation results w.r.t. all the arguments (see, e.g., \cite[Section~I.14]{HNW08}).

\begin{theorem}\label{thm:pertres}
With reference to the solution (\ref{eq:ysol}) of problem (\ref{eq:odexi}), one has:
$$
a)~\frac{\partial}{\partial t}y(t) =f(y(t)),\qquad
b)~\frac{\partial}{\partial \eta} y(t) = \Phi(t,\xi),\qquad
c)~\frac{\partial}{\partial \xi} y(t) = -\Phi(t,\xi)f(\eta),
$$ 
where $\Phi(t,\xi)$ is the solution of the variational problem
$$\dot\Phi(t,\xi) = F(y(t))\Phi(t,\xi), \qquad t\in[\xi,T], \qquad \Phi(\xi,\xi)= I\in\RR^{m\times m},$$     
having set%\rnote{\blue eq:F}
\begin{equation}\label{eq:F}
F(y) = \frac{\partial}{\partial y}f(y).
\end{equation} 
\end{theorem}

From this theorem, the following result readily follows where, hereafter, $|\cdot|$ will denote any convenient vector norm.

\begin{corollary}\label{cor:pert} With reference to (\ref{eq:odexi}), and assuming that $\xi\in[t_{n-1},t_n]$, one has:
$$y(t,\xi,\eta+\delta\eta)= y(t,\xi,\eta) + \Phi(t,\xi)\delta\eta + (t-\xi)O(|\delta\eta|^2), \qquad t\in[t_{n-1},t_n].$$
\end{corollary} 

\subsection{Main results (ODE case)}\label{sec:mainode} With reference to (\ref{eq:h})--(\ref{eq:yn-1}), we are now in the position of 
stating the results concerning the approximation error at the grid points,%\rnote{\blue eq:en}
\begin{equation}\label{eq:en}
y(t_n)-y_n\equiv \hat\sigma_n(h)-\sigma_n(h),\qquad n=1,\dots,N,
\end{equation}
and, more in general, on each subinterval $[t_{n-1},t_n]$:%\rnote{\blue eq:deltasig}
\begin{equation}\label{eq:deltasig}
\delta\sigma_n(ch) := \hat\sigma_n(ch)-\sigma_n(ch) \equiv \hat\sigma(t_{n-1}+ch)-\sigma(t_{n-1}+ch), \qquad c\in[0,1].
\end{equation}

For the first step of the approximation procedure, the following theorem holds true, the proof being similar to that of \cite[Theorem~1]{BIT2012}.

\begin{theorem}\label{thm:ode1}
With reference to  (\ref{eq:en}) and (\ref{eq:deltasig}), one has:
$$y(t_1) - y_1 = O(h^{2s+1}), \qquad \|\delta\sigma_1\| := \max_{c\in[0,1]} |\hat\sigma_1(ch)-\sigma_1(ch)| = O(h^{s+1}).$$
\end{theorem}
\proof  By virtue of Corollary~\ref{cor:gamj} and Theorem~\ref{thm:pertres}, one has:
\begin{eqnarray*}\lefteqn{
\hat\sigma_1(ch)-\sigma_1(ch)~=~y(t_0+ch,t_0,y_0)-y(t_0+ch,t_0+ch,\sigma_1(ch))}\\
 &=& y(t_0+ch,t_0,\sigma_1(0))-y(t_0+ch,t_0+ch,\sigma_1(ch))=\int_{ch}^0 \frac{\dd}{\dd t} y(t_0+ch,t_0+t,\sigma_1(t))\dd t \\
 &=& \int_{ch}^0 \left[ \left.\frac{\partial}{\partial\xi}y(t_0+ch,\xi,\sigma_1(t))\right|_{\xi=t_0+t}
+\left.\frac{\partial}{\partial\eta}y(t_0+ch,t_0+t,\eta)\right|_{\eta=\sigma_1(t)}\dot\sigma_1(t)\right]\dd t\\
&=& \int_0^{ch} \Phi(t_0+ch,t_0+t)\left[f(\sigma_1(t)) - \dot\sigma_1(t) \right]\dd t \\ 
&=& h \int_0^c \Phi(t_0+ch,t_0+\zeta h)\left[f(\sigma_1(\zeta h)) - \dot\sigma_1(\zeta h) \right]\dd \zeta\\
&=&h\int_0^c\Phi(t_0+ch,t_0+\zeta h)\left[ \sum_{j\ge0} P_j(\zeta) \gamma_j(\sigma_1) - \sum_{j=0}^{s-1} P_j(\zeta) \gamma_j(\sigma_1)\right]\dd\zeta\\
&=&h\sum_{j\ge s} \left[\int_0^cP_j(\zeta)\Phi(t_0+ch,t_0+\zeta h)\dd\zeta\right]\underbrace{\gamma_j(\sigma_1)}_{=\,O(h^j)}.
\end{eqnarray*}
Consequently, the second part of the statement follows for $c\in(0,1)$, whereas, when $c=1$ one deduces, by virtue of Theorem~Theorem~\ref{thm:Gh}:
\begin{eqnarray*}
y(t_1)-y_1 &\equiv& \hat\sigma_1(h)-\sigma_1(h) \,=\, h\sum_{j\ge s} \overbrace{\underbrace{\left[\int_0^1P_j(\zeta)\underbrace{\Phi(t_1,t_0+\zeta h)}_{=:\,G(\zeta h)}\dd\zeta\right]}_{=O(h^j)}\gamma_j(\sigma_1)}^{=\,O(h^{2j})}
~=~%\\ &=& 
O(h^{2s+1}).\,\QED
\end{eqnarray*}

For the remaining steps, the following result holds true.

\begin{theorem}\label{thm:oden}
With reference to  (\ref{eq:en}) and (\ref{eq:deltasig}), for $n=1,\dots,N$ one has:
$$y(t_n) - y_n = y(t_{n-1})-y_{n-1}+O(h^{2s+1}), \quad \|\delta\sigma_n\| := \max_{c\in[0,1]} |\delta\sigma_n(ch)| = O(h^{s+1}).$$
\end{theorem}
\proof  By induction on $n$. For $n=1$ the statement follows from the previous Theorem~\ref{thm:ode1}. Assuming it true for $n-1$, for $n$ one has:
\begin{eqnarray*}
\hat\sigma_n(ch) - \sigma_n(ch) &=& y(t_{n-1}+ch,t_{n-1},\hat\sigma_n(0)) - y(t_{n-1}+ch,t_{n-1}+ch,\sigma_n(ch)) \\[2mm]
&=& \underbrace{y(t_{n-1}+ch,t_{n-1},\hat\sigma_n(0)) - y(t_{n-1}+ch,t_{n-1},\sigma_n(0))}_{=:\,E_{n,1}(ch)} ~+\\[1mm]
&&   \underbrace{y(t_{n-1}+ch,t_{n-1},\sigma_n(0)) - y(t_{n-1}+ch,t_{n-1}+ch,\sigma_n(ch))}_{=:\,E_{n,2}(ch)}.
\end{eqnarray*}
By using similar arguments as those used in the proof of Theorem~\ref{thm:ode1}, one deduces that
$$E_{n,2}(ch) = \left\{\begin{array}{cc} O(h^{s+1}), &c\in(0,1),\\[2mm] O(h^{2s+1}), &c=1.\end{array}\right.$$
Moreover, considering that, by the induction hypothesis,
$$\delta\sigma_n(0) = \hat\sigma_n(0)-\sigma_n(0) = y(t_{n-1})-y_{n-1} = (n-1)\,O(h^{2s+1}),$$
from Corollary~\ref{cor:pert}, one has:
\begin{eqnarray*}
E_{n,1}(ch) &=& \underbrace{\Phi(t_{n-1}+ch,t_{n-1})}_{=\,I+O(ch)}\delta\sigma_n(0) + ch\,O(|\delta\sigma_n(0)|^2) \\
&=& y(t_{n-1})-y_{n-1} + (n-1)O(c h^{2s+2}). 
\end{eqnarray*}
Consequently, for $c=1$ one obtains the first part of the statement, whereas the second part follows by taking $c\in(0,1)$.\,\QED

\begin{remark}\label{rem:hbvminf} 
We observe that the two equivalent equations (see (\ref{eq:gammaj}), (\ref{eq:sig1}), and (\ref{eq:sig})):
$$\dot\sigma_n(ch) = \sum_{j=0}^{s-1} P_j(c) \int_0^1 P_j(\zeta) f(\sigma_n(\zeta h))\dd\zeta, \qquad c\in[0,1], \qquad \sigma_n(0) = y_{n-1},$$
and\,%\rnote{\blue eq:MFE}
\begin{equation}\label{eq:MFE}
\sigma_n(ch) = y_{n-1} +h\sum_{j=0}^{s-1} \int_0^cP_j(x)\dd x \int_0^1 P_j(\zeta) f(\sigma_n(\zeta h))\dd\zeta, \qquad c\in[0,1],
\end{equation}
define a so called HBVM$(\infty,s)$ method on the interval $[t_{n-1},t_n]$ (equation (\ref{eq:MFE}) is named {\em Master Functional Equation} in \cite{BIT2010}. See also \cite{BI2016,BI2018}). Consequently, such method defines an order $2s$ approximation procedure for all $s\ge1$, which can be also recast as a continuous-stage Runge-Kutta method \cite{ABI2019}. 
In particular, the case $s=1$ corresponds to the so called AVF method \cite{QMcL2008}; the case $s\ge1$ has been also considered in \cite{Ha2010}. 
\end{remark}

An interesting question concerns the difference between the Fourier coefficients of the solution (\ref{eq:hsig})--(\ref{eq:gammaj}) and those of the polynomial approximation (\ref{eq:sig}) on the interval $[t_{n-1},t_n]$. The next result clarifies the issue.

\begin{theorem}\label{thm:diffgam}
With reference to  (\ref{eq:hsig}), (\ref{eq:gammaj}), and (\ref{eq:sig}), for all $n=1,\dots,N$ one has:
$$\delta\gamma_j^n := \gamma_j(\hat\sigma_n)-\gamma_j(\sigma_n) = O(h^{2s-j}), \qquad j=0,\dots,s-1.$$
\end{theorem}

\proof  First of all, from (\ref{eq:leg}), (\ref{eq:hsig}), (\ref{eq:sig}), and Theorem~\ref{thm:oden} we know that:
$$y(t_n)-y_n = y(t_{n-1})-y_{n-1} + h[\gamma_0(\hat\sigma_n)-\gamma_0(\sigma_n)] = y(t_{n-1})-y_{n-1} + O(h^{2s+1}).$$
Consequently, from the last equality one derives:
$$\delta\gamma_0^n = \gamma_0(\hat\sigma_n)-\gamma_0(\sigma_n) = O(h^{2s}).$$
Further, by taking into account (\ref{eq:F}) and (\ref{eq:deltasig}), one obtains:
\begin{eqnarray*}
O(h^{2s})&=&\gamma_0(\hat\sigma_n)-\gamma_0(\sigma_n) ~=~ \int_0^1 \left[f(\hat\sigma_n(\zeta h))-f(\sigma_n(\zeta h))\right]\dd\zeta\\
 &=& \int_0^1 \underbrace{\int_0^1 F(\sigma_n(\zeta h)+ c\,\delta\sigma_n(\zeta h))\dd c}_{=:\,G(\zeta h)}\delta\sigma_n(\zeta h)\dd\zeta ~=~ \int_0^1 G(\zeta h) \delta\sigma_n(\zeta h)\dd\zeta.
\end{eqnarray*}
Now, considering that ~$P_0(x)\equiv 1$~ and, for all $\zeta\in[0,1]$,%\rnote{\blue eq:legint}
\begin{eqnarray}\label{eq:legint}
\int_0^\zeta P_j(x)\dd x &=& \xi_{j+1}P_{j+1}(\zeta)-\xi_j P_{j-1}(\zeta), \qquad j\ge1,\\
 &&with\qquad \xi_j = \left(2\sqrt{4j^2-1}\right)^{-1},\nonumber
\end{eqnarray}
one has:
\begin{eqnarray*} 
\delta\sigma_n(\zeta h) &=& \hat\sigma_n(\zeta h)-\sigma_n(\zeta h) \\&=& y(t_{n-1})- y_{n-1} \,+\, h\sum_{j= 0}^{s-1}\int_0^\zeta P_j(x)\dd x\,\delta\gamma_j^n \,+\, h\sum_{j\ge s}\int_0^\zeta P_j(x)\dd x\,\gamma_j(\hat\sigma_n)\\
&=& y(t_{n-1})- y_{n-1}\,+\, \zeta h \delta\gamma_0^n + h\sum_{j=1}^{s-1} \left[\xi_{j+1}P_{j+1}(\zeta)-\xi_jP_{j-1}(\zeta)\right]\delta\gamma_j^n\,+\\
&&h\sum_{j\ge s} \left[\xi_{j+1}P_{j+1}(\zeta)-\xi_jP_{j-1}(\zeta)\right]\gamma_j(\hat\sigma_n).  
\end{eqnarray*}
Consequently, from Theorem~\ref{thm:Gh} and Corollary~\ref{cor:gamj},  one obtains:
\begin{eqnarray*}\lefteqn{
O(h^{2s}) ~=~\int_0^1 G(\zeta h)\delta\sigma_n(\zeta h)\,\dd\zeta~=~\underbrace{\int_0^1 G(\zeta h)\dd\zeta}_{=\,O(1)}\underbrace{[y(t_{n-1})- y_{n-1}]}_{=\,(n-1)\,O(h^{2s+1})}}\\
&&+~h\underbrace{\int_0^1  G(\zeta h)\zeta\dd\zeta}_{=\,O(1)}\, \underbrace{\delta\gamma_0^n}_{=\,O(h^{2s})}~+
~h \sum_{j=1}^{s-1} \int_0^1  G(\zeta h)\left[\xi_{j+1}P_{j+1}(\zeta)-\xi_jP_{j-1}(\zeta)\right]\dd\zeta\,\delta\gamma_j^n \\
&&+\underbrace{h\sum_{j\ge s} \overbrace{\int_0^1  G(\zeta h)\left[\xi_{j+1}P_{j+1}(\zeta)-\xi_jP_{j-1}(\zeta)\right]\dd\zeta}^{=\,O(h^{j-1})}\,\overbrace{\gamma_j(\hat\sigma_n)}^{=\,O(h^j)}}_{=\,O(h^{2s})},
\end{eqnarray*}
from which, 
$$O(h^{2s})\,=\,h \sum_{j=1}^{s-1} \overbrace{\int_0^1  G(\zeta h)\left[\xi_{j+1}P_{j+1}(\zeta)-\xi_jP_{j-1}(\zeta)\right]\dd\zeta}^{=\,O(h^{j-1})}\delta\gamma_j^n$$
follows and, therefore, one concludes that ~$\delta\gamma_j^n=O(h^{2s-j})$, ~$j=0\dots,s-1$.\,\QED

\subsection{Conservative/dissipative problems}\label{sec:codipro}
An interesting case \cite{FM2010,BI2016,MB2016} is that when problem (\ref{eq:ode1}) is in the form%\nnote{\blue eq:codipro}
\begin{equation}\label{eq:codipro}
\dot y(t) = S\nabla H(y(t)), \qquad t\in[t_0,T], \qquad y(t_0)=y_0\in\RR^m,
\end{equation}
with $S\in\RR^{m\times m}$ either a skew-symmetric matrix, $S^\top=-S$, or a negative semidefinite matrix, $S\le0$, whereas $\nabla H$ is the gradient of a scalar function usually called the {\em Hamiltonian}. As is clear:
\begin{itemize}
\item when $S^\top=-S$: $$\frac{\dd}{\dd t} H(y(t)) = \nabla H(y(t))^\top \dot y(t) =  \nabla H(y(t))^\top S \nabla H(y(t)) = 0,$$
so that $H$ is a conserved quantity, and the problem is said to be {\em conservative};

\item when $S\le 0$: $$\frac{\dd}{\dd t} H(y(t)) = \nabla H(y(t))^\top \dot y(t) =  \nabla H(y(t))^\top S \nabla H(y(t)) \le 0,$$
and the problem is said to be {\em dissipative}.
\end{itemize}
The next result shows that this behavior is preserved by the approximations (\ref{eq:sig1})--(\ref{eq:yn-1}), upon observing that in this case (\ref{eq:gammaj}) can be conveniently rewritten as%\nnote{\blue eq:gamjcodipro}
\begin{equation}\label{eq:gamjcodipro}
\gamma_j(z) = S  \int_0^1 P_j(\zeta) \nabla H(z(\zeta h))\dd\zeta =: S\beta_j(z).
\end{equation}

\begin{theorem}\label{thm:codipro}
With reference to (\ref{eq:sig1})--(\ref{eq:yn-1}) applied for approximating problem (\ref{eq:codipro}), for all $n=1,\dots,N$ one has:
\begin{itemize}
\smallskip
\item $H(y_n) = H(y_{n-1})$, ~when~ $S^\top = -S$; 

\smallskip
\item $H(y_n) \le H(y_{n-1})$, ~when~ $S\le0$.
\end{itemize}
\end{theorem}
\proof  In fact, by considering that $y_n=\sigma_n(h)$, $y_{n-1}=\sigma_n(0)$, and taking into account (\ref{eq:gamjcodipro}), one has: 
\begin{eqnarray*}\lefteqn{
H(y_n)-H(y_{n-1})~=~ H(\sigma_n(h))-H(\sigma_n(0)) = \int_0^h \frac{\dd}{\dd t}H(\sigma_n(t))\dd t} \\
&=& h\int_0^1 \nabla H(\sigma_n(c h))^\top\dot\sigma_n(c h)\dd c =  h\int_0^1 \nabla H(\sigma_n(c h))^\top\sum_{j=0}^{s-1} P_j(c)\gamma_j(\sigma_n)\dd c\\
&=&  h\sum_{j=0}^{s-1}\underbrace{\left[\int_0^1 \nabla H(\sigma_n(c h)) P_j(c)\dd c\right]^\top}_{=\,\beta_j(\sigma_n)^\top} S\beta_j(\sigma_n) =
h\sum_{j=0}^{s-1}\beta_j(\sigma_n)^\top S\beta_j(\sigma_n)~=:~\Delta H_n.
\end{eqnarray*}
Consequently, if $S^\top=-S$, then $\Delta H_n=0$, whereas $\Delta H_n\le0$, when $S\le0$.\,\QED

\begin{remark}\label{rem:hbvminfcodipro}
According to Remark~\ref{rem:hbvminf}, one then obtains that HBVM$(\infty,s)$ methods can preserve the conservative/dissipative feature of problem (\ref{eq:codipro}).
\end{remark}

\subsection{Discretization and Runge-Kutta formulation}\label{sec:odedis}
Quoting Dahlquist and Bj\"ork \cite[p.\,521]{DB2008} {\em``as is well known, even many relatively simple integrals cannot be expressed in finite terms of elementary functions, and thus must be evaluated by numerical methods.''} In this context, this quite obvious statement means that the approximation procedure defined by (\ref{eq:sig1}) and (\ref{eq:gammaj}) does not yet provide a ``true'' numerical method. In fact, the integrals defining the Fourier coefficients,%\rnote{\blue eq:fourier}
\begin{equation}\label{eq:fourier}
\gamma_j(\sigma_n) = \int_0^1 P_j(\zeta)f(\sigma_n(\zeta h))\dd\zeta, \qquad j=0,\dots,s-1,
\end{equation} 
need to be numerically approximated by using a quadrature formula. Since we are dealing with a polynomial approximation, it is quite natural to do this by using an interpolatory quadrature with abscissae and weights $(c_i,b_i)$, $i=1,\dots,k$ (we shall always assume $k$ distinct abscissae):%\rnote{\blue eq:quad}
\begin{equation}\label{eq:quad}
\gamma_j(\sigma_n) = \sum_{i=1}^k b_i P_j(c_i)f(\sigma_n(c_i h)) + \Delta_j(h), 
\end{equation}
where $\Delta_j(h)$ is the quadrature error. The following straightforward result holds true.

\begin{theorem}\label{thm:quaderr} If the quadrature $(c_i,b_i)$, $i=1,\dots,k$ has order $q$, i.e., it is exact for polynomial integrands of degree $q-1$, then $$\Delta_j(h) = O(h^{q-j}), \qquad j=0,\dots,s-1.$$
\end{theorem}

\begin{remark}
As is well known, since the quadrature (\ref{eq:quad}) is based at $k$ (distinct) abscissae, $q\in\{k,\dots,2k\}$: the lower limit is obtained by a generic choice of the abscissae, whereas the upper one is achieved by placing them at the zeros of $P_k(c)$.
\end{remark}

When using a quadrature, clearly the Fourier coefficients (\ref{eq:fourier}) may be not exactly evaluated anymore. This implies that we are actually computing a possibly different piecewise polynomial approximation $u(t)$ such that (compare with (\ref{eq:gammaj})--(\ref{eq:yn-1})), for all $n=1,\dots,N$:
%\nnote{\blue eq:un eq:un1}
\begin{eqnarray}\label{eq:un}
u_n(ch) &\equiv& u(t_{n-1}+ch),\qquad c\in[0,1],\\ \label{eq:un1}
\dot u_n(ch) &=& \sum_{j=0}^{s-1} P_j(c) \hat\gamma_j(u_n), \qquad c\in[0,1],\qquad u_n(0) = y_{n-1},
\end{eqnarray}
with (see (\ref{eq:quad}))%\nnote{\blue eq:hgamj-unc-yn-2}
\begin{eqnarray}\label{eq:hgamj}
\hat\gamma_j(u_n) &:=& \sum_{i=1}^k b_i P_j(c_i) f(u_n(c_i h))\equiv \gamma_j(u_n)-\Delta_j(h), \quad j=0,\dots,s-1,\\ \label{eq:unc}
u_n(ch) &=& y_{n-1} + h\sum_{j=0}^{s-1} \int_0^c P_j(x)\dd x\, \hat\gamma_j(u_n), \qquad c\in[0,1],\\
u_n(h) &=:& y_n \equiv y_{n-1} + h\hat\gamma_0(u_n).\label{eq:yn-2}
\end{eqnarray}
Actually, (\ref{eq:hgamj})--(\ref{eq:yn-2}) define the $n$th integration step, by using a timestep $h$, performed with the $k$ stage Runge-Kutta method having stages:%\nnote{\blue eq:Yin}
\begin{equation}\label{eq:Yin}
Y_i^n := u_n(c_ih), \qquad i=1,\dots,k.
\end{equation} 
In fact, evaluating (\ref{eq:unc}) at the abscissae $c_1,\dots,c_k$, and substituting in it the $s$ approximate Fourier coefficients (\ref{eq:hgamj}), one obtains, after rearranging terms,%\nnote{\blue eq:but1 eq:but2}
\begin{eqnarray}\label{eq:but1}
Y_i^n &=& y_{n-1} + h\sum_{\ell=1}^k \underbrace{b_\ell \sum_{j=0}^{s-1} \int_0^{c_i} P_j(x)\dd x\, P_j(c_\ell)}_{=:\,a_{i\ell}} f(Y_\ell^n), \qquad i=1,\dots,n,\\ \label{eq:but2}
y_n &=& y_{n-1} + h\sum_{i=1}^k b_i f(Y_i^n).
\end{eqnarray}
In other words, we have derived the $k$-stage Runge-Kutta method with abscissae and weights $(c_i,b_i)$, $i=1,\dots,k$, and Butcher matrix
$A=\begin{pmatrix} a_{i\ell}\end{pmatrix}\in\RR^{k\times k}$. Next theorem puts the Butcher tableau in a more compact form \cite{BI2016}.

\begin{theorem}\label{thm:but} The Butcher tableau of the Runge-Kutta method (\ref{eq:but1})--(\ref{eq:but2}) is given by%\nnote{\blue eq:but3}
\begin{equation}\label{eq:but3}
\begin{array}{c|c}
\bfc & \I_s\P_s^\top\Omega\\
\hline
&\\[-3.5mm]
&\bfb^\top 
\end{array}
\end{equation}
where
$$\bfb = \begin{pmatrix} b_1&\dots&b_k\end{pmatrix}^\top, \qquad 
\bfc = \begin{pmatrix} c_1&\dots&c_k\end{pmatrix}^\top, \qquad \Omega = \diag(\bfb),$$
and
$$\P_s = \begin{pmatrix} P_0(c_1)&\dots&P_{s-1}(c_1)\\ \vdots & &\vdots \\ P_0(c_k)&\dots&P_{s-1}(c_k)\end{pmatrix},\qquad
\I_s = \begin{pmatrix} \int_0^{c_1}P_0(x)\dd x&\dots&\int_0^{c_1}P_{s-1}(x)\dd x\\ \vdots & &\vdots \\\int_0^{c_k}P_0(x)\dd x&\dots&\int_0^{c_k}P_{s-1}(x)\dd x\end{pmatrix}.$$
\end{theorem}

It is possible to derive an alternative formulation of the Runge-Kutta method (\ref{eq:but3}).
In fact, using the relation (\ref{eq:legint}) between the integrals of the Legendre polynomials and the polynomials themselves, and considering that
$$\int_0^c P_0(x)\dd x =  \xi_1P_1(c) + \xi_0P_0(c), \qquad \xi_0=\frac{1}2,$$ one has that \,$\I_s = \P_{s+1}\hat X_s$, where
$$\P_{s+1} = \begin{pmatrix} P_0(c_1)&\dots&P_s(c_1)\\ \vdots & &\vdots \\ P_0(c_k)&\dots&P_s(c_k)\end{pmatrix}, 
\qquad
\hat X_s = \begin{pmatrix} \xi_0 &-\xi_1\\
                                   \xi_1  &0  &\ddots\\
                                             &\ddots &\ddots & -\xi_{s-1}\\
                                                         &           &\xi_{s-1}  &0\\
                                                         \hline
                                                         &           &            &\xi_s \end{pmatrix} =: \begin{pmatrix} X_s\\ \hline 0\,\dots\,0\,\xi_s\end{pmatrix}.$$
Consequently, the Butcher tableau (\ref{eq:but3}) can be rewritten as 
$$ 
\begin{array}{c|c}
\bfc & \P_{s+1}\hat X_s\P_s^\top\Omega\\
\hline
&\\[-3.5mm]
&\bfb^\top 
\end{array}~.
$$ 
When the quadrature (\ref{eq:quad}) has order $q\ge 2s$, it is quite straightforward to prove that
$$\P_s^\top\Omega \P_s = I_s, \qquad \P_s^\top\Omega \P_{s+1} = [ I_s ~ \bfzero]\in\RR^{s\times (s+1)},$$
where in general, hereafter, $I_r\in\RR^{r\times r}$ is the identity matrix (when the dimension of the identity matrix is not explicitly indicated, it will be easily deducible from the context). Consequently,
$$\P_s^\top\Omega\left[  \P_{s+1}\hat X_s\P_s^\top\Omega \right] \P_s= X_s,$$
which can be regarded as a generalization of the $W$-transformation in \cite[Theorem\,5.6, p.\,79]{HW2002}. 
In addition to this, when $q\ge 2s$ also the following results hold true (for sake of brevity, we do not discuss the case $q<2s$, since it has no practical interest).

\begin{theorem}\label{thm:ode1_1}
With reference to (\ref{eq:ode1})--(\ref{eq:gammaj}) and  (\ref{eq:un})--(\ref{eq:yn-2}), 
and assuming that the quadrature formula (\ref{eq:quad}) has order $q\ge 2s$, one has:
$$y(t_1) - y_1 = O(h^{2s+1}), \qquad  \max_{c\in[0,1]} |\hat\sigma_1(ch)-u_1(ch)| = O(h^{s+1}).$$
\end{theorem}

\begin{theorem}\label{thm:oden_1}
With reference to (\ref{eq:ode1})--(\ref{eq:gammaj}) and  (\ref{eq:un})--(\ref{eq:yn-2}), 
and assuming that the quadrature formula (\ref{eq:quad}) has order $q\ge 2s$, for $n=1,\dots,N$ one has:
$$y(t_n) - y_n = y(t_{n-1})-y_{n-1}+O(h^{2s+1}), \qquad \max_{c\in[0,1]} |\hat\sigma_n(ch)-u_n(ch)| = O(h^{s+1}).$$
\end{theorem}

\begin{theorem}\label{thm:diffgam_1}
With reference to (\ref{eq:hsig})--(\ref{eq:gammaj}) and (\ref{eq:hgamj})--(\ref{eq:yn-2}),  and assuming that the quadrature formula (\ref{eq:quad}) has order $q\ge 2s$, for all $n=1,\dots,N$ one has:
$$\delta\hat\gamma_j^n:=\gamma_j(\hat\sigma_n)-\hat\gamma_j(u_n) = O(h^{2s-j}), \qquad j=0,\dots,s-1.$$
\end{theorem}

Concerning the case of conservative/dissipative problems in the form (\ref{eq:codipro}), the result of Theorem~\ref{thm:codipro} modifies as follows.

\begin{theorem}\label{thm:codipro_1}
With reference to (\ref{eq:un})--(\ref{eq:yn-2}) applied for approximating problem (\ref{eq:codipro}), and assuming that the quadrature formula (\ref{eq:quad}) has order $q\ge 2s$, for all $n=1,\dots,N$ one has:
\begin{itemize}
\smallskip
\item if $H$ is a polynomial of degree not larger than $q/s$, then the result of Theorem~\ref{thm:codipro} continues to hold; 

\smallskip
\item differently,\begin{itemize}
\item $H(y_n) = H(y_{n-1})+O(h^{q+1})$, ~when~ $S^\top = -S$,

\smallskip
\item $H(y_n) \le H(y_{n-1})+O(h^{q+1})$, ~when~ $S\le0$.
\end{itemize}\end{itemize}
\end{theorem}

We here provide only the proof of Theorem~\ref{thm:ode1_1} (see also \cite[Theorem\,4]{BIT2012}), since those of Theorem~\ref{thm:oden_1}, Theorem~\ref{thm:diffgam_1}, and Theorem~\ref{thm:codipro_1} can be similarly obtained by slightly adapting the corresponding proofs of Theorem~\ref{thm:oden}, Theorem~\ref{thm:diffgam}, and Theorem~\ref{thm:codipro}, respectively.

\proof (of Theorem~\ref{thm:ode1_1})
By taking into account the result of Theorem~\ref{thm:quaderr}, one has:
\begin{eqnarray*}\lefteqn{
\hat\sigma_1(ch)-u_1(ch)=y(t_0+ch,t_0,y_0)-y(t_0+ch,t_0+ch,u_1(ch)) }\\ 
&=& y(t_0+ch,t_0,u_1(0))-y(t_0+ch,t_0+ch,u_1(ch))
= \int_{ch}^0 \frac{\dd}{\dd t} y(t_0+ch,t_0+t, u_1(t))\dd t \\
&=& \int_{ch}^0 \left[ \left.\frac{\partial}{\partial\xi}y(t_0+ch,\xi, u_1(t))\right|_{\xi=t_0+t}
+\left.\frac{\partial}{\partial\eta}y(t_0+ch,t_0+t,\eta)\right|_{\eta=u_1(t)}\dot u_1(t)\right]\dd t\\
&=& \int_0^{ch} \Phi(t_0+ch,t_0+t)\left[f(u_1(t)) - \dot u_1(t) \right]\dd t \\ 
&=& h \int_0^c \Phi(t_0+ch,t_0+\zeta h)\left[f(u_1(\zeta h)) - \dot u_1(\zeta h) \right]\dd \zeta\\
&=&h\int_0^c\Phi(t_0+ch,t_0+\zeta h)\left[ \sum_{j\ge0} P_j(\zeta) \gamma_j(u_1) - \sum_{j=0}^{s-1} P_j(\zeta) \hat\gamma_j(u_1)\right]\dd\zeta\\
&=&h\int_0^c\Phi(t_0+ch,t_0+\zeta h)\left[ \sum_{j\ge0} P_j(\zeta) \gamma_j(u_1) - \sum_{j=0}^{s-1} P_j(\zeta) \left(\gamma_j(u_1)-\Delta_j(h)\right)\right]\dd\zeta\\
&=&h\sum_{j\ge s} \left[\int_0^cP_j(\zeta)\Phi(t_0+ch,t_0+\zeta h)\dd\zeta\right]\underbrace{\gamma_j(u_1)}_{=\,O(h^j)} \\
&&\,+\,h\sum_{j=0}^{s-1}\left[\int_0^cP_j(\zeta)\Phi(t_0+ch,t_0+\zeta h)\dd\zeta\right]\underbrace{\Delta_j(h)}_{=\,O(h^{q-j})}.
\end{eqnarray*}
Consequently, the second part of the statement follows by considering that, for $c\in(0,1)$, this quantity is
$$O(h^{s+1}) + O(h^{q-s+2}) = O(h^{s+1}),$$ since $q\ge 2s$, whereas, when $c=1$ one deduces, by virtue of Theorem~\ref{thm:Gh}, and considering that $t_1=t_0+h$:
\begin{eqnarray*}%\lefteqn{
y(t_1)-y_1 &\equiv& \hat\sigma_1(h)-u_1(h) \,=\, h\sum_{j\ge s} \overbrace{\underbrace{\left[\int_0^1P_j(\zeta)\underbrace{\Phi(t_1,t_0+\zeta h)}_{=:\,G(\zeta h)}\dd\zeta\right]}_{=\,O(h^j)}\gamma_j(u_1)}^{=\,O(h^{2j})} 
%} \\ 
\end{eqnarray*}\begin{eqnarray*}%%%%%%%%%%%%%%
&&\,+\,
h\sum_{j=0}^{s-1} \overbrace{\underbrace{\left[\int_0^1P_j(\zeta)\underbrace{\Phi(t_1,t_0+\zeta h)}_{=:\,G(\zeta h)}\dd\zeta\right]}_{=\,O(h^j)}\Delta_j(h)}^{=\,O(h^q)} 
\,=\,O(h^{2s+1})+O(h^{q+1}) ~=~ O(h^{2s+1}).\,\QED
\end{eqnarray*}

\begin{remark}\label{rem:hbvm}
When the $k$ abscissae are placed at the zeros of $P_k(c)$, and $k\ge s$, one obtains a HBVM$(k,s)$ method, whose order is $2s$ \cite{BIT2010,BI2016,BI2018}. It is worth mentioning that the HBVM$(s,s)$ method is nothing but the $s$-stage Gauss-Legendre collocation  method. Moreover, the HBVM$(k,1)$ methods correspond to the second-order Runge-Kutta methods described in \cite{CMcLMcLOQW2009}. Different choices of the quadrature have been also considered in \cite{IP2007,IP2008,IT2009,BIT2015}.
\end{remark}

\subsection{Solving the discrete problems}\label{sec:soldis} Sometimes, the number of stages $k$ of the Runge-Kutta method (\ref{eq:but3}) can be much larger than the degree $s$ of the underlying polynomial approximation (\ref{eq:hgamj})--(\ref{eq:yn-2}). This is the case, for example, of HBVM$(k,s)$ methods when used as energy-conserving methods \cite{BIT2010,BI2016,BI2018} (see also Theorem~\ref{thm:codipro_1} in Section~\ref{sec:codipro}). In such a case, it is clear that the usual implementation of the Runge-Kutta method leads to the solution of a discrete problem having (block) dimension $k$. Nevertheless, for sake of completeness we now recall how the discrete problem to be solved can be actually recast so as to have (block) dimension $s$, independently of $k$ \cite{BIT2011}. This clearly allows for relatively large values of $k$, thus making possible the use of an arbitrarily high-order quadrature (\ref{eq:quad}). Let us then consider the first integration step of the method for solving (\ref{eq:ode1}) with timestep $h$,  (thus, we can skip the index $n$ of the step). Setting $\bfuno=\begin{pmatrix}1\,,\,\dots\,,\,1\end{pmatrix}^\top\in\RR^k$, and $Y$ the stage vector of (block) dimension $k$, one obtains that the stage equation for (\ref{eq:but3}) is given by:%\rnote{\blue eq:stageq}
\begin{equation}\label{eq:stageq}
Y = \bfuno\otimes y_0 + h\I_s\P_s^\top\Omega\otimes I_m\, f(Y),
\end{equation}
with an obvious meaning of $f(Y)$. However, we observe that \cite{BIT2011}
$$ \P_s^\top\Omega\otimes I_m\,f(Y) =: \hat\bfgamma \equiv 
\begin{pmatrix} \hat\gamma_0(u_1)\\ \vdots \\ \hat\gamma_{s-1}(u_1)\end{pmatrix},$$
i.e., the (block) vector with the $s$ coefficients of the polynomial approximation $u_1(ch)$ (see (\ref{eq:hgamj})--(\ref{eq:unc})). Consequently, (\ref{eq:stageq}) can be rewritten as
$$Y = \bfuno\otimes y_0 + h\I_s\otimes I_m\, \hat\bfgamma.$$
By combining the last two equations one eventually obtains:%\rnote{\blue eq:gammeq}
\begin{equation}\label{eq:gammeq}
\hat\bfgamma = \P_s^\top\Omega\otimes I_m\,f\left(\bfuno\otimes y_0 + h\I_s\otimes I_m\, \hat\bfgamma\right),
\end{equation}
which is a discrete problem, equivalent to (\ref{eq:stageq}), having (block) dimension $s$, {\em independently} of $k$. Once this equation has been solved, the new approximation is derived, according to (\ref{eq:yn-2}), as $$y_1 = y_0 + h\hat\gamma_0(u_1).$$
It is also worth mentioning that very effective nonlinear iterations have been devised for solving (\ref{eq:gammeq}) \cite{BIT2011,BFCI2014,BI2016} (the most effective being that derived from the so called {\em blended iteration} introduced in \cite{BM2002}, see also \cite{BM2009}).

\section{The DDE case}\label{sec:dde}
As for the ODE case, also for DDEs we shall consider, without loss of generality, the simpler problem%\rnote{\blue eq:dde1}
\begin{eqnarray}\label{eq:dde1}
\dot y(t) &=& f(y(t),y(t-\tau)), \qquad t\in[t_0,T], \qquad y(t_0) = y_0,\\ \nonumber
y(t)&=&\phi(t), \qquad\qquad\qquad\quad\, t\in[t_0-\tau,t_0),\nonumber
\end{eqnarray} 
in place of (\ref{eq:ddeivp}) where, usually, $y_0=\phi(t_0)$. Moreover, we shall suppose that both the timestep $h$ defining the discrete mesh (\ref{eq:h}) and the width of the integration interval, $T-t_0$, are commensurable with the delay:%\rnote{\blue eq:hnu}
\begin{equation}\label{eq:hnu}
\tau = \nu h, \qquad T-t_0 = K\tau, \qquad K,\nu\in\NN, 
\end{equation}
so that the discrete mesh is now given by:%\rnote{\blue eq:h1}
\begin{equation}\label{eq:h1}
t_n = t_0+nh, \qquad n=-\nu,\dots,N\equiv K\nu.
\end{equation} 

On one hand, similarly as done in the ODE case, let us denote, for notational purposes, by $\hat\sigma(t)\equiv y(t)$ the solution of (\ref{eq:dde1}), and%\rnote{\blue eq:sign-d} 
\begin{equation}\label{eq:sign-d}
\hat\sigma_n(ch) :=\hat\sigma(t_{n-1}+ch), \qquad c\in[0,1], \qquad n=1-\nu,\dots,N,
\end{equation}
its restriction to the time interval $[t_{n-1},t_n]$. Consequently,%\nnote{\blue eq:hsigfi}
\begin{equation}\label{eq:hsigfi}
\hat\sigma_n(ch) \equiv \phi(t_{n-1}+ch), \qquad c\in[0,1], \qquad n=1-\nu,\dots,0,
\end{equation}
whereas, for $n=1,\dots,N$, one has (compare with (\ref{eq:hsig1})--(\ref{eq:gammaj})):%\nnote{\blue eq:sig1-d}
\begin{equation}\label{eq:hsig1-d}
\dot{\hat\sigma}_n(ch) = \sum_{j\ge 0} P_j(c)\gamma_j(\hat\sigma_n,\hat\sigma_{n-\nu}), \qquad c\in[0,1], \qquad \hat\sigma_n(0) = y(t_{n-1}),
\end{equation}
so that,%\nnote{\blue eq:hsig-d}
\begin{equation}\label{eq:hsig-d}
\hat\sigma_n(ch) = y(t_{n-1})+h\sum_{j\ge0} \int_0^c P_j(x)\dd x\,\gamma_j(\hat\sigma_n,\hat\sigma_{n-\nu}), \qquad c\in[0,1],
\end{equation}
and%\nnote{\blue eq:ytn-d}
\begin{equation}\label{eq:ytn-d}
y(t_n) = y(t_{n-1}) + h\gamma_0(\hat\sigma_n,\hat\sigma_{n-\nu}) \equiv \hat\sigma(t_n),
\end{equation}
where, in general, for any suitably regular functions $z,w:[0,h]\rightarrow\RR^m$,%\nnote{\blue eq:gammaj-d}
\begin{equation}\label{eq:gammaj-d}
\gamma_j(z,w) := \int_0^1 P_j(\zeta)f(z(\zeta h),w(\zeta h))\dd\zeta.
\end{equation}

On the other hand, we shall look for a piecewise  approximation to $\hat\sigma(t)$, i.e. $\sigma(t)$, such that (compare with (\ref{eq:sign})--(\ref{eq:yn-1}))%\nnote{\blue eq:sign-1-d}
\begin{equation}\label{eq:sign-1-d}
\sigma_n(ch) := \sigma(t_{n-1}+ch), \qquad c\in[0,1], \qquad n=1-\nu,\dots,N,
\end{equation}
denotes its restriction to the time interval $[t_{n-1},t_n]$. Consequently, one has:%\nnote{\blue eq:sigfi}
\begin{equation}\label{eq:sigfi}
\sigma_n(ch) \equiv \phi(t_{n-1}+ch), \qquad c\in[0,1], \qquad n=1-\nu,\dots,0,
\end{equation}
whereas, for $n=1,\dots,N$, $\sigma_n\in\Pi_s$ satisfies the differential equation%\nnote{\blue eq:sig1-d}
\begin{equation}\label{eq:sig1-d}
\dot\sigma_n(ch) = \sum_{j=0}^{s-1} P_j(c)\gamma_j(\sigma_n,\sigma_{n-\nu}), \qquad c\in[0,1], \qquad \sigma_n(0) = y_{n-1},
\end{equation}
so that,%\nnote{\blue eq:sig-d}
\begin{equation}\label{eq:sig-d}
\sigma_n(ch) = y_{n-1}+h\sum_{j=0}^{s-1} \int_0^c P_j(x)\dd x\,\gamma_j(\sigma_n,\sigma_{n-\nu}), \qquad c\in[0,1],
\end{equation}
and%\nnote{\blue eq:ytn-1-d}
\begin{equation}\label{eq:ytn-1-d}
y_n = y_{n-1} + h\gamma_0(\sigma_n,\sigma_{n-\nu}) =:\sigma_n(h),
\end{equation}
with $\gamma_j(\sigma_n,\sigma_{n-\nu})$ defined according to (\ref{eq:gammaj-d}). In the sequel, we shall discuss the accuracy of the approximations:%\nnote{\blue eq:appro-d}
\begin{eqnarray}\label{eq:appro-d}
y(t_n)-y_n&\equiv& \hat\sigma_n(h)-\sigma_n(h), \\[1mm] \nonumber
 \delta\sigma_n(ch) &:=&\hat\sigma_n(ch)-\sigma_n(ch), \quad c\in(0,1),\qquad n=1,\dots,N.
\end{eqnarray}
For this purpose, some preliminary results are given in the next section.

\subsection{Preliminary results}\label{sec:preldde} We start with the generalization of Corollary~\ref{cor:gamj} to the present setting.

\begin{corollary}\label{cor:gamj-d} With reference to (\ref{eq:gammaj-d}), one has: \,$\gamma_j(z,w) = O(h^j)$.
\end{corollary}
\proof 
Immediate from Theorem~\ref{thm:Gh}, by setting $G(\zeta h):=f(z(\zeta h),w(\zeta h))$.\,\QED\bigskip

We also need perturbation results corresponding to those of Theorem~\ref{thm:pertres} for ODEs. For this purpose, it is sufficient to discuss them for a local problem defined on two contiguous time subintervals of width $\tau$: the former containing the {\em memory}, the latter containing the solution to be computed. Without loss of generality, we shall then fix the reference interval $[t_0-\tau,t_0+\tau]$, where we consider the following problem, defined for a generic $\xi\in [t_0,t_0+\tau]$:%\rnote{\blue eq:dde1loc}
\begin{eqnarray}\label{eq:dde1loc}
\dot y(t) &=& f(y(t),y(t-\tau)), \qquad t\in[t_0,t_0+\tau], \qquad y(\xi)=\eta\in\RR^m,\\
y(t)&=&\phi(t), \qquad\qquad\qquad\quad\, t\in[t_0-\tau,t_0).\nonumber
\end{eqnarray}
Problem (\ref{eq:dde1loc}) defines a generalization of the localized one associated to (\ref{eq:dde1}) (obtained for $\xi=t_0$ and $\eta=y_0$), and we shall denote its solution by\,
%\footnote{As done in the ODE case, for sake of brevity, hereafter we shall use either one of the two notations, depending on the needs.}%\rnote{\blue eq:ytpar}
\begin{equation}\label{eq:ytpar}
y(t)\equiv y(t,\xi,\eta,\phi;t_0),
\end{equation}
in order to emphasize its dependence on the first four parameters, whereas the last one refers to the time subinterval.
We shall also use the following notation:%\rnote{\blue eq:F12}
\begin{equation}\label{eq:F12}
F_1(z,w) = \frac{\partial}{\partial z}f(z,w), \qquad F_2(z,w) = \frac{\partial}{\partial w}f(z,w).
\end{equation}

\begin{remark}\label{rem:memory}
It is clear that the function $\phi$ in (\ref{eq:dde1loc}) represents the {\em memory term} of the equation, and it is a known function. The same will happen in the subsequent reference interval, $[t_0,t_0+2\tau]$, obtained by shifting to the right the previous one by $\tau$, once the solution of (\ref{eq:dde1loc}) has been computed, and so forth.
\end{remark}

To begin with, let us state the following straightforward result, whose proof is omitted for brevity.

\begin{theorem}\label{thm:tutto}
The solution (\ref{eq:ytpar}) of problem (\ref{eq:dde1loc}) is defined on the whole time interval $[t_0,t_0+\tau]$, independently of the point $\xi\in[t_0,t_0+\tau]$ where the condition $\eta$ is given. 
\end{theorem}

The following result then holds true (compare with Theorem~\ref{thm:pertres}).

\begin{theorem}\label{thm:pertres-d}
With reference to the solution (\ref{eq:ytpar}) of problem (\ref{eq:dde1loc}), one has: 
$$
a)~\frac{\partial}{\partial t}y(t) =f(y(t),y(t-\tau)),\quad
b)~\frac{\partial}{\partial \eta} y(t) = \Phi(t,\xi;t_0), \quad %$$ $$
c)~\frac{\partial}{\partial \xi} y(t) = -\Phi(t,\xi;t_0)f(y(\xi),y(\xi-\tau)),
$$ 
where $\Phi(t,\xi;t_0)$ satisfies (see (\ref{eq:F12})):%\rnote{\blue eq:varpro1}
\begin{eqnarray}\nonumber
%\dot\Phi(t,\xi;t_0) &=& F_1(y(t),y(t-\tau))\Phi(t,\xi;t_0)+F_2(y(t),y(t-\tau))\Phi(t-\tau,\xi;t_0), \\
%&&\qquad\qquad\qquad\qquad\qquad\qquad\qquad\qquad\qquad\qquad\qquad t\in[t_0,t_0+\tau], \\
\dot\Phi(t,\xi;t_0) &=& F_1(y(t),y(t-\tau))\Phi(t,\xi;t_0), \qquad t\in[t_0,t_0+\tau], \\[1mm] \label{eq:varpro1}
\Phi(\xi,\xi;t_0) &=&I\in\RR^{m\times m}, \\[1mm] \nonumber
\Phi(t,\xi;t_0)&=&O\in\RR^{m\times m}, \qquad t\in[t_0-\tau,t_0).
\end{eqnarray}   
\end{theorem}
\proof 
The statement $a)$ clearly follows from (\ref{eq:dde1loc}). From the same equation one also derives that, for $t\in[t_0,t_0+\tau]$, %\nnote{\blue eq:F2}
\begin{eqnarray}\label{eq:F2}
 \qquad \frac{\dd}{\dd t}\left(\frac{\partial}{\partial\eta}y(t)\right) &=& \frac{\partial}{\partial\eta} \dot y(t)~=~ \frac{\partial}{\partial\eta} f(y(t),y(t-\tau))\\
 &=& F_1(y(t),y(t-\tau))\frac{\partial}{\partial\eta} y(t)\,+\,F_2(y(t),y(t-\tau))\frac{\partial}{\partial\eta} y(t-\tau).\nonumber
\end{eqnarray}
Moreover, at $t=\xi$, 
$$\frac{\partial}{\partial\eta} y(\xi) = \frac{\partial}{\partial\eta} \eta = I,$$
and, for $t\in[t_0-\tau,t_0)$,
$$\frac{\partial}{\partial\eta} y(t) = \frac{\partial}{\partial\eta} \phi(t) = O.$$
This latter equality implies that, for $t\in[t_0,t_0+\tau)$, the term $F_2(y(t),y(t-\tau))\frac{\partial}{\partial\eta} y(t-\tau)$ in (\ref{eq:F2}) vanishes, thus reducing to the first equation in (\ref{eq:varpro1}), so that $b)$ eventually follows. Finally,  by virtue of Theorem~\ref{thm:tutto}, let $t^*$ be a generic point in the interval $[t_0,t_0+\tau]$, and denote $$y^*=y(t^*,\xi,\eta,\phi;t_0).$$ Consequently, since $\xi\in[t_0,t_0+\tau]$ as well, one has: 
$$\eta = y(\xi,t^*,y(t^*,\xi,\eta,\phi;t_0),\phi;t_0),$$ 
so that we eventually arrive at the identity
$$y^* = y(t^*,\xi,y(\xi,t^*,y^*,\phi;t_0),\phi;t_0).$$
By taking into account the results of the previous points $a)$ and $b)$, one derives:
\begin{eqnarray*}
0 &=& \frac{\dd}{\dd \xi}y(t^*,\xi,y(\xi,t^*,y^*,\phi;t_0),\phi;t_0)\\
&=& \frac{\partial}{\partial \xi}y(t^*,\xi,\eta,\phi;t_0) + \frac{\partial}{\partial \eta}y(t^*,\xi,\eta,\phi;t_0) \left.\frac{\partial}{\partial t}y(t,t^*,y^*,\phi;t_0)\right|_{t=\xi}\\
&=& \frac{\partial}{\partial \xi}y(t^*,\xi,\eta,\phi;t_0) + \Phi(t^*,\xi;t_0) f(y(\xi),y(\xi-\tau)).
\end{eqnarray*}
The statement $c)$ then follows, by taking into account that $t^*$ is generic.\,\QED\bigskip

One main difference with the ODE case, stems from the fact that now (\ref{eq:ytpar}) also depends on the {\em memory term} $\phi$, which is a {\em functional} parameter. Consequently, we now look for a Frech\'et derivative  such that, for any perturbation
 $\delta\phi\in C([t_0-\tau,t_0])$ and $t\in[t_0-\tau,t_0+\tau]$: %\rnote{\red referenze? io ho messo [20, App.\,A]}%\nnote{\blue eq:freder}
 \begin{equation}\label{eq:freder}
\lim_{\eps\rightarrow0}\frac{y(t,\xi,\eta,\phi+\eps\delta\phi;t_0) - y(t,\xi,\eta,\phi;t_0)}{\eps}= \int_{t_0-\tau}^{t_0} \frac{\delta}{\delta \phi(\zeta)}y(t)\delta \phi(\zeta)\dd\zeta,
\end{equation}
where%\nnote{\blue eq:funder} 
\begin{equation}\label{eq:funder}
\frac{\delta}{\delta \phi(\zeta)}y(t) \,:\,(t,\zeta)\in[t_0-\tau,t_0+\tau]\times [t_0-\tau,t_0)\rightarrow  \begin{pmatrix}  \frac{\delta}{\delta \phi_j(\zeta)}y_i(t) \end{pmatrix}\in\RR^{m\times m},
\end{equation}
is the {\em functional derivative} of (\ref{eq:ytpar}) (see, e.g., \cite[Appendix\,A]{ED2011}), with $y_i$ and $\phi_j$ the respective entries of $y$ and $\phi$.  For later use, we recall that, for a given $\hat t\in[t_0-\tau,t_0)$ and $i,j=1,\dots,m$,%\rnote{\blue eq:funderij}
\begin{eqnarray}\label{eq:funderij}
~~~\frac{\delta}{\delta \phi_j(\hat t\,)}y_i(t) &\equiv& \lim_{\eps\rightarrow0} \frac{y_i(t,\xi,\eta,\phi+\eps e_j\delta_{\hat t};t_0) - y_i(t,\xi,\eta,\phi;t_0)}{\eps} \\
%&=& \int_{t_0-\tau}^{t_0} \left[\frac{\delta}{\delta \phi(\zeta)}y_i(t)\right]^\top e_j\delta_{\hat t}(\zeta)\dd\zeta
&=&\int_{t_0-\tau}^{t_0} \frac{\delta}{\delta \phi_j(\zeta)}y_i(t)\delta_{\hat t}(\zeta)\dd\zeta,  \nonumber
\end{eqnarray}
with $e_j\in\RR^m$ the $j$th unit vector and, hereafter, $\delta_{\hat t}(t)$ is the Dirac delta function centered at $\hat t$. The following result holds true.

\begin{lemma}\label{lem:xit} With reference to (\ref{eq:freder}) and (\ref{eq:funder}), for any fixed $t\in[\xi,t_0+\tau]\subseteq[t_0,t_0+\tau]$ one has:
$$\frac{\delta}{\delta\phi(\zeta)}y(t)=O\in\RR^{m\times m}, \qquad \forall\zeta\in[t_0-\tau,\xi-\tau)\cup(t-\tau,t_0).$$
\end{lemma} 
\proof  Having fixed $t\in[\xi,t_0+\tau]$, it follows that $\forall\zeta\in[t_0-\tau,\xi-\tau)\cup(t-\tau,t_0)$, setting as usual $\delta_\zeta$ the Dirac delta centered at $\zeta$, one has:
$$y(t,\xi,\eta,\phi+\eps\delta_\zeta; t_0) = y(t,\xi,\eta,\phi; t_0), \qquad \forall\eps\in\RR.$$
In fact, by virtue of (\ref{eq:dde1loc}), the solution (\ref{eq:ytpar}) is independent of the values of $\phi$ outside the interval $[\xi-\tau,t-\tau]$. Consequently, by  taking into account (\ref{eq:funderij}), it follows that:
$$\frac{\delta}{\delta\phi(\zeta)}y(t)=\lim_{\eps\rightarrow0} \frac{y(t,\xi,\eta,\phi+\eps\delta_\zeta; t_0) - y(t,\xi,\eta,\phi; t_0)}{\eps}=O.\,\QED$$
  
\begin{figure}[t]
\centering
 \centerline{ \includegraphics[width=1\textwidth,keepaspectratio]{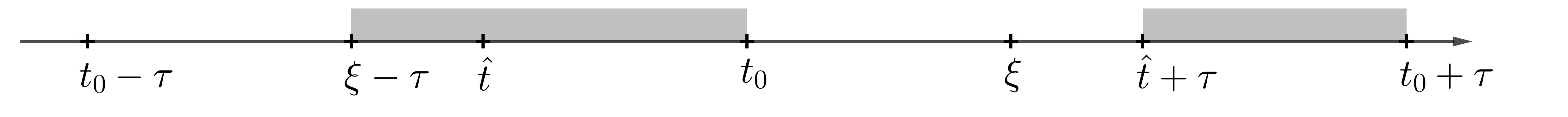}}
      \vspace{-3mm}
  \caption{Relevant time subintervals for Theorem~\ref{thm:frechet}.}
  \label{fig:figura}
\end{figure}

Taking into account Lemma~\ref{lem:xit}, the following result provides a more practical characterization of the functional derivative (\ref{eq:freder})--(\ref{eq:funder}).
 Figure~\ref{fig:figura} displays the location of the most relevant points and subintervals involved in Theorem~\ref{thm:frechet}.
% which may be more easily understood by making reference to \ref{fig:figura}.

\begin{theorem}\label{thm:frechet}
With reference to the solution (\ref{eq:ytpar}) of problem (\ref{eq:dde1loc}), for any $\hat t\in(\xi-\tau,t_0)$ one has:
\begin{equation}\label{eq:vdpsi}
\frac{\delta}{\delta \phi(\hat t\,)} y(t) = \Psi(t,\hat t\,;t_0),
\end{equation}   
where $\Psi(t,\hat t\,;t_0)$ satisfies (see (\ref{eq:F12})):%\rnote{\blue eq:varpro2}
\begin{eqnarray}\nonumber
 \dot\Psi(t,\hat t\,;t_0) &=& F_1(y(t),y(t-\tau))\Psi(t,\hat t\,;t_0),\qquad t\in(\hat t+\tau,t_0+\tau],\\[1mm]  \label{eq:varpro2} 
\quad\Psi(\hat t+\tau,\hat t\,;t_0) &=& F_2(y(\hat t+\tau),y(\hat t\,)),\\ [1mm]
\Psi(t,\hat t;t_0)&=&\delta_{\hat t}(t)I, \qquad t\in[t_0-\tau,\hat t+\tau),\nonumber
\end{eqnarray}
with $O,I\in\RR^{m\times m}$ and  \,$\delta_{\hat t}(t)$\, the Dirac delta function.
\end{theorem}
\proof  In fact, for $t\in[t_0,t_0+\tau]$ one has, by virtue of (\ref{eq:dde1loc}):
\begin{eqnarray*}\lefteqn{
\frac{\dd}{\dd t}\frac{\delta}{\delta \phi(\hat t\,)} y(t,\xi,\eta,\phi;t_0) ~=~\frac{\delta}{\delta \phi(\hat t\,)} \dot y(t)  ~=~ \frac{\delta}{\delta \phi(\hat t\,)} f(y(t),y(t-\tau))} \\
&=& F_1(y(t),y(t-\tau)) \frac{\delta}{\delta \phi(\hat t\,)} y(t) + F_2(y(t),y(t-\tau)) \frac{\delta}{\delta \phi(\hat t\,)} y(t-\tau),
\end{eqnarray*}
i.e., using the notation (\ref{eq:vdpsi}),
\begin{eqnarray}\label{eq:vdpsi0}
\qquad \dot\Psi(t,\hat t\,;t_0) &=& F_1(y(t),y(t-\tau))\Psi(t,\hat t\,;t_0) + F_2(y(t),y(t-\tau))\Psi(t-\tau,\hat t\,;t_0).\\[-3mm] \nonumber
\end{eqnarray} 
Moreover, at $t=\xi$,
\begin{eqnarray}\label{eq:psixi0}
\Psi(\xi,\hat t\,;t_0) &=& \frac{\delta}{\delta \phi(\hat t\,)} y(\xi) =  \frac{\delta}{\delta \phi(\hat t\,)} \eta = O,\\[-3mm] \nonumber
\end{eqnarray}
since the condition $y(\xi)=\eta$ is independent of the history $\phi$. Further, taking into account (\ref{eq:funderij}), for all $i,j=1,\dots,m$ and  for $t\in[t_0-\tau,t_0)$, one has:
$$ \frac{\delta}{\delta \phi_j(\hat t\,)} y_i(t) =  \frac{\delta}{\delta \phi_j(\hat t\,)} \phi_i(t) = \lim_{\eps\rightarrow0}\frac{[\phi_i(t)+\eps \delta_{ij}\delta_{\hat t}(t)]-\phi_i(t)}{\eps} = \delta_{ij}\delta_{\hat t}(t),$$
with $\delta_{ij}$ the Kronecker delta. Consequently,
\begin{eqnarray}\label{eq:psit0}
\Psi(t,\hat t\,;t_0) &=& \frac{\delta}{\delta \phi(\hat t\,)} y(t) = \delta_{\hat t}(t)I, \qquad t\in[t_0-\tau,t_0). \\[-3mm] \nonumber
\end{eqnarray}
From (\ref{eq:vdpsi0})--(\ref{eq:psit0}), one then derives, considering that (see Figure~\ref{fig:figura}) \,$t_0-\tau\le\xi-\tau<\hat t:$
%\begin{equation}\label{eq:psiht} 
\begin{eqnarray}\label{eq:psiht} 
\quad \Psi(t,\hat t\,;t_0) &=& \int_\xi^t \dot \Psi(\zeta,\hat t\,;t_0)\dd\zeta = \left\{\begin{array}{cc} O, & t\in[t_0,\hat t+\tau),\\[2mm]
F_2(y(\hat t+\tau),y(\hat t\,)), &t=\hat t+\tau.\end{array}\right. \\[-3mm] \nonumber
\end{eqnarray}
%\end{equation}
From (\ref{eq:psit0}) and (\ref{eq:psiht}) the last two equations in (\ref{eq:varpro2}) follow. Consequently, from (\ref{eq:vdpsi0}), one obtains
$$\dot\Psi(t,\hat t\,;t_0) = F_1(y(t),y(t-\tau))\Psi(t,\hat t\,;t_0), \qquad t\in(\hat t+\tau,t_0+\tau],$$
which completes the proof of (\ref{eq:varpro2}).\,\QED\bigskip

As a straightforward consequence, the following result holds true, which guarantees the regularity of $\Psi$ w.r.t. its first two arguments (again, for sake of clarity, refer to Figure~\ref{fig:figura}).

\begin{corollary}\label{cor:troppofico}
With reference to the solution (\ref{eq:ytpar}) of problem (\ref{eq:dde1loc}), and considering (\ref{eq:F12}), (\ref{eq:varpro1}), and (\ref{eq:varpro2}), for any $\hat t\in(\xi-\tau,t_0)$ one has:
\begin{equation}\label{eq:psifi}
\Psi(t,\hat t\,;t_0) = \Phi(t,\hat t+\tau; t_0)F_2(y(\hat t+\tau),y(\hat t\,)), \qquad t\in[\hat t+\tau,t_0+\tau].
\end{equation}
\end{corollary}

Finally, the following result holds true (compare with Corollary~\ref{cor:pert} of the ODE case).

\begin{corollary}\label{cor:pert-d} With reference to the solution (\ref{eq:ytpar}) of problem (\ref{eq:dde1loc}), and considering (\ref{eq:varpro1}) and (\ref{eq:psifi}), for any $\delta\phi\in C([t_0-\tau,t_0])$ one has: 
\begin{eqnarray*}
y(t,\xi,\eta+\delta\eta,\phi+\delta\phi;t_0) &=& y(t,\xi,\eta,\phi;t_0) + \Phi(t,\xi;t_0)\delta\eta+\int_{\xi-\tau}^{t-\tau} \Psi(t,\zeta;t_0)\delta\phi(\zeta)\dd\zeta \\ &&+\,(t-\xi)\,O(|\delta\eta|+\|\delta\phi\|)^2, \qquad\qquad  t\in[\xi,t_0+\tau],
\end{eqnarray*}
with \,$\|\delta\phi\|=\max_{\zeta\in[\xi-\tau,t-\tau]}|\delta\phi(\zeta)|$. 
\end{corollary}
\proof 
The statement follows from Theorem~\ref{thm:pertres-d}, part $b)$, and Theorem~\ref{thm:frechet}, by taking into account (\ref{eq:freder}) and the result of Lemma~\ref{lem:xit}.\,\QED\bigskip

\subsection{Main results (DDE case)}\label{sec:maindde} We are now in the position of discussing the accuracy of the approximations (\ref{eq:appro-d}). To begin with, the following result holds true.

\begin{theorem}\label{thm:ddem}
With reference to  (\ref{eq:dde1})--(\ref{eq:appro-d}), for $n=1,\dots,\nu$ one has:
$$y(t_n) - y_n = y(t_{n-1})-y_{n-1}+O(h^{2s+1}), \quad \|\delta\sigma_n\| := \max_{c\in[0,1]} |\delta\sigma_n(ch)| = O(h^{s+1}).$$
\end{theorem}
\proof  The statement follows from Theorem~\ref{thm:oden} by considering that, for $n=1,\dots,\nu$, $t\in[t_0,t_0+\tau]$ in (\ref{eq:dde1}), so that $y(t-\tau)\equiv \phi(t-\tau)$, which is a known function, thus obtaining an ODE.\,\QED\bigskip

This result allows us to state the following one, which generalizes that of Theorem~\ref{thm:diffgam} to the present case.

\begin{theorem}\label{thm:diffgam-d}
With reference to  (\ref{eq:hsig-d}), (\ref{eq:gammaj-d}), (\ref{eq:sig-d}), and (\ref{eq:appro-d}) if for $n\ge 1$ one has:
$$y(t_r) - y_r = y(t_{r-1})-y_{r-1}+O(h^{2s+1}), \qquad r=1,\dots,n,$$
then
$$\delta\gamma_j^n := \gamma_j(\hat\sigma_n,\hat\sigma_{n-\nu})-\gamma_j(\sigma_n,\sigma_{n-\nu}) = O(h^{2s-j}), \qquad j=0,\dots,s-1.$$
\end{theorem}
\proof  The proof is by generalized induction. 
For $n=1,\dots,\nu$\, the statement follows from Theorem~\ref{thm:diffgam} and Theorem~\ref{thm:ddem} since, in this case, 
$$\hat\sigma_{n-\nu}(ch)\equiv\sigma_{n-\nu}(ch)\equiv \phi(t_{n-1}+ch-\tau), \qquad c\in[0,1],$$
so that we are dealing with an ODE. Assume now it true up to $n-1$, and prove for $n$. By hypothesis, and from (\ref{eq:ytn-d}) and (\ref{eq:ytn-1-d}), we know that
$$y(t_n) - y_n = y(t_{n-1})-y_{n-1}+O(h^{2s+1}) \equiv y(t_{n-1})-y_{n-1}+h\delta\gamma_0^n,$$ so that $\delta\gamma_0^n = O(h^{2s})$ follows. Then, by taking into account (\ref{eq:F12}), it follows that:
\begin{eqnarray*}
O(h^{2s})&=& \delta\gamma_0^n \,=\,\gamma_0(\hat\sigma_n,\hat\sigma_{n-\nu})-\gamma_0(\sigma_n,\sigma_{n-\nu}) \\
&=& \int_0^1 \left[f(\hat\sigma_n(\zeta h),\hat\sigma_{n-\nu}(\zeta h))-f(\sigma_n(\zeta h),\sigma_{n-\nu}(\zeta h))\right]\dd\zeta\\
 &=&\int_0^1 \int_0^1 \left[\, F_1(\sigma_n(\zeta h)+ c\,\delta\sigma_n(\zeta h),\sigma_{n-\nu}(\zeta h)+ c\,\delta\sigma_{n-\nu}(\zeta h))\delta\sigma_n(\zeta h) ~+\right.\\[1mm]
&&\qquad\quad~ \left.F_2(\sigma_n(\zeta h)+ c\,\delta\sigma_n(\zeta h),\sigma_{n-\nu}(\zeta h)+ c\,\delta\sigma_{n-\nu}(\zeta h))\delta\sigma_{n-\nu}(\zeta h) \,\right]\dd c\,\dd\zeta \\[2mm]
&=&\int_0^1 \underbrace{\int_0^1 F_1(\sigma_n(\zeta h)+ c\,\delta\sigma_n(\zeta h),\sigma_{n-\nu}(\zeta h)+ c\,\delta\sigma_{n-\nu}(\zeta h))\dd c}_{=:\,G_1(\zeta h)}\, \delta\sigma_n(\zeta h)\,\dd \zeta ~+\\
&&\int_0^1 \underbrace{\int_0^1 F_2(\sigma_n(\zeta h)+ c\,\delta\sigma_n(\zeta h),\sigma_{n-\nu}(\zeta h)+ c\,\delta\sigma_{n-\nu}(\zeta h))\dd c}_{=:\,G_2(\zeta h)}\,\delta\sigma_{n-\nu}(\zeta h)\,\dd \zeta \\
&=& \int_0^1 G_1(\zeta h)\delta\sigma_n(\zeta h)\dd\zeta \,+\,  \int_0^1 G_2(\zeta h)\delta\sigma_{n-\nu}(\zeta h)\dd\zeta.
\end{eqnarray*}
Let us discuss in detail the term
$$\int_0^1 G_1(\zeta h)\delta\sigma_n(\zeta h)\,\dd \zeta,$$
since the remaining one,
$$\int_0^1 G_2(\zeta h)\delta\sigma_{n-\nu}(\zeta h)\,\dd \zeta = O(h^{2s}),$$
is similarly discussed, by taking into account the induction hypothesis. By virtue of (\ref{eq:legint}), one has:
\begin{eqnarray*} 
\delta\sigma_n(\zeta h) &=& \hat\sigma_n(\zeta h)-\sigma_n(\zeta h) \\&=& y(t_{n-1})- y_{n-1} \,+\, h\sum_{j= 0}^{s-1}\int_0^\zeta P_j(x)\dd x\,\delta\gamma_j^n \,+\, \sum_{j\ge s}\int_0^\zeta P_j(x)\dd x\,\gamma_j(\hat\sigma_n,\hat\sigma_{n-\nu})\\
&=& y(t_{n-1})- y_{n-1}\,+\, \zeta h \delta\gamma_0^n + h\sum_{j=1}^{s-1} \left[\xi_{j+1}P_{j+1}(\zeta)-\xi_jP_{j-1}(\zeta)\right]\delta\gamma_j^n\,+\\
&&h\sum_{j\ge s} \left[\xi_{j+1}P_{j+1}(\zeta)-\xi_jP_{j-1}(\zeta)\right]\gamma_j(\hat\sigma_n,\hat\sigma_{n-\nu}).  
\end{eqnarray*}
Consequently, from Theorem~\ref{thm:Gh} and Corollary~\ref{cor:gamj-d},  one obtains:
\begin{eqnarray*}\lefteqn{
O(h^{2s}) ~=~\int_0^1 G_1(\zeta h)\delta\sigma_n(\zeta h)\,\dd\zeta~=~\underbrace{\int_0^1 G_1(\zeta h)\dd\zeta}_{=\,O(1)}\underbrace{[y(t_{n-1})- y_{n-1}]}_{=\,(n-1)\,O(h^{2s+1})}}\\
&&+\,h\underbrace{\int_0^1  G_1(\zeta h)\zeta\dd\zeta}_{=\,O(1)}\, \underbrace{\delta\gamma_0^n}_{=\,O(h^{2s})}\,+\,h \sum_{j=1}^{s-1} \int_0^1  G_1(\zeta h)\left[\xi_{j+1}P_{j+1}(\zeta)-\xi_jP_{j-1}(\zeta)\right]\dd\zeta\,\delta\gamma_j^n \\
&&+\underbrace{h\sum_{j\ge s} \overbrace{\int_0^1  G_1(\zeta h)\left[\xi_{j+1}P_{j+1}(\zeta)-\xi_jP_{j-1}(\zeta)\right]\dd\zeta}^{=\,O(h^{j-1})}\,\overbrace{\gamma_j(\hat\sigma_n,\hat\sigma_{n-\nu})}^{=\,O(h^j)}}_{=\,O(h^{2s})},
\end{eqnarray*}
from which, 
$$O(h^{2s}) = h \sum_{j=1}^{s-1} \underbrace{\int_0^1  G_1(\zeta h)\left[\xi_{j+1}P_{j+1}(\zeta)-\xi_jP_{j-1}(\zeta)\right]\dd\zeta}_{=\,O(h^{j-1})}\delta\gamma_j^n$$
follows and, therefore, one concludes that ~$\delta\gamma_j^n=O(h^{2s-j})$, ~$j=0\dots,s-1$.\,\QED\bigskip

As a consequence, the following result can be stated.

\begin{theorem}\label{thm:ddeN}
With reference to  (\ref{eq:dde1})--(\ref{eq:appro-d}), for $n=1,\dots,N\equiv K\nu$ one has:
$$y(t_n) - y_n = y(t_{n-1})-y_{n-1}+O(h^{2s+1}), \quad \|\delta\sigma_n\| := \max_{c\in[0,1]} |\delta\sigma_n(ch)| = O(h^{s+1}).$$
\end{theorem}
\proof  The proof is done by induction on groups of $\nu$ consecutive steps. For the first $\nu$ steps, the statement follows from Theorem~\ref{thm:ddem}.
Assume now, by induction, that it holds true up to $t_{k\nu}=t_0+k\nu h$, and let us prove for $n=k\nu+1,\dots,(k+1)\nu$. For this purpose, for $k=1,\dots,K-1$\, let us set:
$$\phi_k(t) \equiv \sigma(t), \qquad \hat\phi_k(t)\equiv \hat\sigma(t)\equiv y(t), \qquad t\in[t_{(k-1)\nu},t_{k\nu}).$$
Assuming, again, true the statement for $n-1$, and using the notation (\ref{eq:ytpar}), one has:
\begin{eqnarray*}\lefteqn{\delta\sigma_n(ch)\,=\,\hat\sigma_n(ch)-\sigma_n(ch)}\\[1mm]
&=& y(t_{n-1}+ch,t_{n-1},y(t_{n-1}),\hat\phi_k;t_{k\nu})-y(t_{n-1}+ch,t_{n-1}+ch,\sigma_n(ch),\phi_k;t_{k\nu})\\%[1mm]
&=&\underbrace{y(t_{n-1}+ch,t_{n-1},\overbrace{\sigma_n(0)}^{=\,y_{n-1}},\phi_k;t_{k\nu})-y(t_{n-1}+ch,t_{n-1}+ch,\sigma_n(ch),\phi_k;t_{k\nu})}_{=:\,E_{n,1}^{(k)}(ch)} \\%[1mm] 
&&+~\underbrace{y(t_{n-1}+ch,t_{n-1},y(t_{n-1}),\hat\phi_k;t_{k\nu})-y(t_{n-1}+ch,t_{n-1},\overbrace{y_{n-1}}^{=\,\sigma_n(0)},\phi_k;t_{k\nu})}_{=:\,E_{n,2}^{(k)}(ch)}. 
\end{eqnarray*}
From Theorem~\ref{thm:ddem}, it follows that
%\begin{equation}\label{eq:Ek1}
\begin{eqnarray}\label{eq:Ek1} 
E_{n,1}^{(k)}(h) &=& \sigma_n(0)-\sigma_n(0) + O(h^{2s+1}) \,=\,O(h^{2s+1}),\\[1mm] 
\|E_{n,1}^{(k)}\|&:=&\max_{c\in[0,1]}|E_{n,1}^{(k)}(ch)|=O(h^{s+1}).\nonumber
\end{eqnarray}
%\end{equation}
Moreover, from Corollary~\ref{cor:pert-d}, and considering that $h\nu=\tau$, one has:
\begin{eqnarray*}
E_{n,2}^{(k)}(ch) &=& \overbrace{\Phi(t_{n-1}+ch,t_{n-1};t_{k\nu})}^{=\,I+O(ch)} \overbrace{\left[y(t_{n-1})-y_{n-1}\right]}^{=\,\delta\sigma_n(0)} \\
&&+\, h\int_0^c\Psi(t_{n-1}+ch,t_{n-1}+\zeta h-\tau;t_{k\nu})\delta\sigma_{n-\nu}(\zeta h)\dd\zeta\\
&&+\, ch\, O(|\delta\sigma_n(0)|+\|\delta\sigma_{n-\nu}\|)^2.
\end{eqnarray*}
By considering that $$|\delta\sigma_n(0)|=(n-1)O(h^{2s+1}), \qquad \|\delta\sigma_{n-\nu}\|=O(h^{s+1}),$$
one eventually derives 
$$\|E_{n,2}^{(k)}\|:=\max_{c\in[0,1]}|E_{n,2}^{(k)}(ch)|=O(h^{s+2}),$$
from which the second part of the statement follows, by taking into account (\ref{eq:Ek1}). Moreover, when $c=1$ then $t_{n-1}+h=t_n$ and, by virtue of Theorem~\ref{thm:Gh} and Theorem~\ref{thm:diffgam-d}, one obtains:
\begin{eqnarray*}\lefteqn{
\int_0^1\overbrace{\Psi(t_n,t_{n-1}+\zeta h-\tau;t_{k\nu})}^{=:\,G(\zeta h)}\delta\sigma_{n-\nu}(\zeta h)\dd\zeta}\\
 &=& \int_0^1 G(\zeta h) \left[\sum_{j=0}^{s-1} P_j(\zeta)\delta\gamma_j^{n-\nu}+\sum_{j\ge s} P_j(\zeta) \gamma_j(\hat\sigma_{n-\nu},\hat\sigma_{n-2\nu})\right]\\
&=& \sum_{j=0}^{s-1} \underbrace{\int_0^1P_j(\zeta)G(\zeta h)\dd\zeta}_{=\,O(h^j)} \underbrace{\delta\gamma_j^{n-\nu}}_{=\,O(h^{2s-j})}+\sum_{j\ge s} \underbrace{\int_0^1P_j(\zeta)G(\zeta h)\dd\zeta}_{=\,O(h^j)} \underbrace{\gamma_j(\hat\sigma_{n-\nu},\hat\sigma_{n-2\nu})}_{=\,O(h^j)} %\\&=&
~=~O(h^{2s}).
\end{eqnarray*}
Consequently,
$$E_{n,2}^{(k)}(h) = y(t_{n-1})-y_{n-1} + O(h^{2s+1}),$$
and also the first part of the statement follows.\,\QED 

\subsection{Discretization}\label{ddedis}
The discretization issue proceeds as in the ODE case. In fact, also in the DDE case, the Fourier coefficients (see (\ref{eq:gammaj-d}) and (\ref{eq:sig1-d})),
$$\gamma_j(\sigma_n,\sigma_{n-\nu}) = \int_0^1 P_j(\zeta)f(\sigma_n(\zeta h),\sigma_{n-\nu}(\zeta h))\dd\zeta,\qquad j=0,\dots,s-1,$$
need to be approximated by using a (interpolatory) quadrature rule of order $q$, thus providing a possibly different piecewise approximation $u(t)$, 
$$u(t)\equiv \phi(t), \quad t<t_0, \qquad u_n(ch) := u(t_{n-1}+ch), \quad c\in[0,1],\quad n=1-\nu,\dots,N,$$
such that, for $n\ge 1$:%\nnote{\blue eq:u1-d}
\begin{equation}\label{eq:u1-d}
\dot u_n(ch) = \sum_{j=0}^{s-1} P_j(c)\hat\gamma_j(u_n,u_{n-\nu}), \qquad c\in[0,1], \qquad u_n(0) = y_{n-1}.
\end{equation}
Consequently,%\nnote{\blue eq:u-d}
\begin{equation}\label{eq:u-d}
u_n(ch) = y_{n-1}+h\sum_{j=0}^{s-1} \int_0^c P_j(x)\dd x\,\hat\gamma_j(u_n,u_{n-\nu}), \qquad c\in[0,1],
\end{equation}
and%\nnote{\blue eq:un-1-d}
\begin{equation}\label{eq:un-1-d}
y_n = y_{n-1} + h\hat\gamma_0(u_n,u_{n-\nu}) =: u_n(h),
\end{equation}
where (see (\ref{eq:gammaj-d})),%\nnote{\blue eq:hgammaj-d}
\begin{equation}\label{eq:hgammaj-d}
\hat\gamma_j(u_n,u_{n-\nu}) := \sum_{i=1}^k b_i P_j(c_i)f(u_n(c_ih),u_{n-\nu}(c_ih)) = \gamma_j(u_n,u_{n-\nu}) -\Delta_j(h),
\end{equation}
with $(c_i,b_i)$ the abscissae and weights of the quadrature, and $\Delta_j(h)=O(h^{q-j})$ the quadrature error, where $q$ is the order of the quadrature. 

Formulae (\ref{eq:u1-d}) and (\ref{eq:gammaj-d}) form a subclass of the so called natural continuous RK methods for DDEs (see \cite[Sec. 6.2]{BZ03}). As a consequence, their convergence properties could be as well derived by more classical approaches such as Bellman's method of steps, which is an analytic procedure specific for DDEs. In the present context, the main goal is to show how the framework based on the perturbation theory applied to the truncated Fourier expansion is easily adapted to cope with DDEs, therefore we will pursue this  route of investigation. A further strength of this approach is the possibility of analyzing the convergence properties of the truncated Fourier approximations when these are used as spectral methods in time. In this regard, the analysis for the ODE case has been addressed in \cite{ABI2020}, while a spectral implementation of the methods for DDEs has been considered in \cite{BIZ2021}.

By using standard arguments (which we omit, as done in the ODE case), we can derive the following results, representing the corresponding counterparts of Theorem~\ref{thm:diffgam-d} and Theorem~\ref{thm:ddeN}, respectively.

\begin{theorem}\label{thm:diffgam1-d}
With reference to  (\ref{eq:hsig-d}), (\ref{eq:gammaj-d}), (\ref{eq:u1-d})--(\ref{eq:un-1-d}),
and assuming that the quadrature formula (\ref{eq:hgammaj-d}) has order $q\ge 2s$, if for $n\ge 1$ one has:
$$y(t_r) - y_r = y(t_{r-1})-y_{r-1}+O(h^{2s+1}), \qquad r=1,\dots,n,$$
then
$$\delta\hat\gamma_j^n := \gamma_j(\hat\sigma_n,\hat\sigma_{n-\nu})-\hat\gamma_j(u_n,u_{n-\nu}) = O(h^{2s-j}), \qquad j=0,\dots,s-1.$$
\end{theorem}

\begin{theorem}\label{thm:ddeN-1}
With reference to  (\ref{eq:hsig-d}), (\ref{eq:gammaj-d}), (\ref{eq:u1-d})--(\ref{eq:un-1-d}),
and assuming that the quadrature formula (\ref{eq:hgammaj-d}) has order $q\ge 2s$,  for $n=1,\dots,N\equiv K\nu$ one has:
$$y(t_n) - y_n = y(t_{n-1})-y_{n-1}+O(h^{2s+1}), \qquad \max_{c\in[0,1]} |\hat\sigma_n(ch)-u_n(ch)| = O(h^{s+1}).$$
\end{theorem}

\begin{remark} It is worth mentioning that the result of Theorem~\ref{thm:ddeN-1} states that the super-convergence order $2s$ at the mesh-points $t_n$ is obtained, even though possibly different Runge-Kutta methods are used at each integration step, provided that they define a polynomial approximation of degree $s$. This, in turn, represents a generalization of the results in \cite{B84} for collocation methods.
\end{remark}

We conclude this section,  by recalling that the considerations in Remark~\ref{rem:hbvm} continue to hold in the DDE case and by observing that, concerning the implementation of the resulting Runge-Kutta method used for solving problem (\ref{eq:dde1}), the arguments in Section~\ref{sec:soldis}, {\em mutatis mutandis}, apply as well.

\section{Numerical tests}\label{sec:numtest}
In this section we report a few numerical tests for the DDE case. In fact, in the ODE case, HBVMs have been extensively used as energy-conserving methods for Hamiltonian systems (see, e.g., \cite{ABI2022,BI2016,BI2018,BIMR2019,BMR2019}). We show that, under some circumstances, their use can be advantageous also in the DDE case.  Hereafter, we consider a class of DDEs defined by a Hamiltonian function
\begin{equation}\label{eq:Hdde}
H:(q,p)\in\RR^m\times \RR^m \rightarrow \RR,
\end{equation}
through the equations
\begin{eqnarray}\nonumber
\dot q(t) &=& \quad H_p(q(t),p(t)) +\alpha H_p(q(t-\tau),p(t-\tau)), \\[-2mm] \label{eq:qpdde}\\ \nonumber
 \dot p(t) &=& -\left[H_q(q(t),p(t)) +\alpha H_q(q(t-\tau),p(t-\tau))\right],\quad
\end{eqnarray}
with $\alpha$ a real parameter, $\tau>0$ the delay, and $H_q$ and $H_p$ the partial derivatives of $H$ w.r.t. $q$ and $p$, respectively.  The problem is completed by the initial conditions
\begin{equation}\label{eq:inidde}
q(t) = \phi(t), \qquad p(t) = \psi(t), \qquad t\in[-\tau,0].
\end{equation}
The introduction of such a kind of {\em delay Hamiltonian system} is partly inspired by the problem  of looking for periodic orbits of DDEs, which has been attacked by many authors in the past (see, e.g., \cite{DRB1986,KY1974,Nu1973,Nu1973_1,Nu1979,PaNu2011,W1975}). In this respect, the first two examples below show an attractive  periodic orbit with integer period lying on a level set of the Hamiltonian function (\ref{eq:Hdde}) which is, therefore, a constant of motion once the periodic orbit has been approached. In the third example we are instead interested in simulating the correct qualitative behavior of a {\em dissipative Hamiltonian delay problem} in the phase space when the dynamics takes place in a neighborhood of a separatrix. Taking aside a theoretical discussion of problem  (\ref{eq:Hdde})-(\ref{eq:qpdde}) which would go beyond the scopes of the present work, we infer its properties for the three considered instances by preliminarily applying a high order integrator with very small stepsize, in order to get a very accurate numerical solution that will be taken as a reference trajectory in the phase space. 

For all the three problems, we show that a very accurate approximation of the Hamiltonian function allows us to reproduce the correct geometric features of the solution in the discrete setting.  To the best of our knowledge, this is the first instance of the use of HBVMs in the context of DDEs displaying geometric properties. For comparison purposes, we also solve the problems with the classical Gauss collocation integrator of the same order.
The numerical tests have been implemented in Matlab (R2020b) on a 3 GHz Intel Xeon W10 core computer with 64GB of memory. 

\subsubsection*{Problem~1} % dhbvm.m + Htau.m
\begin{figure}[hp]
\centering

\vspace{-2.5cm}
\centerline{ \includegraphics[width=.55\textwidth]{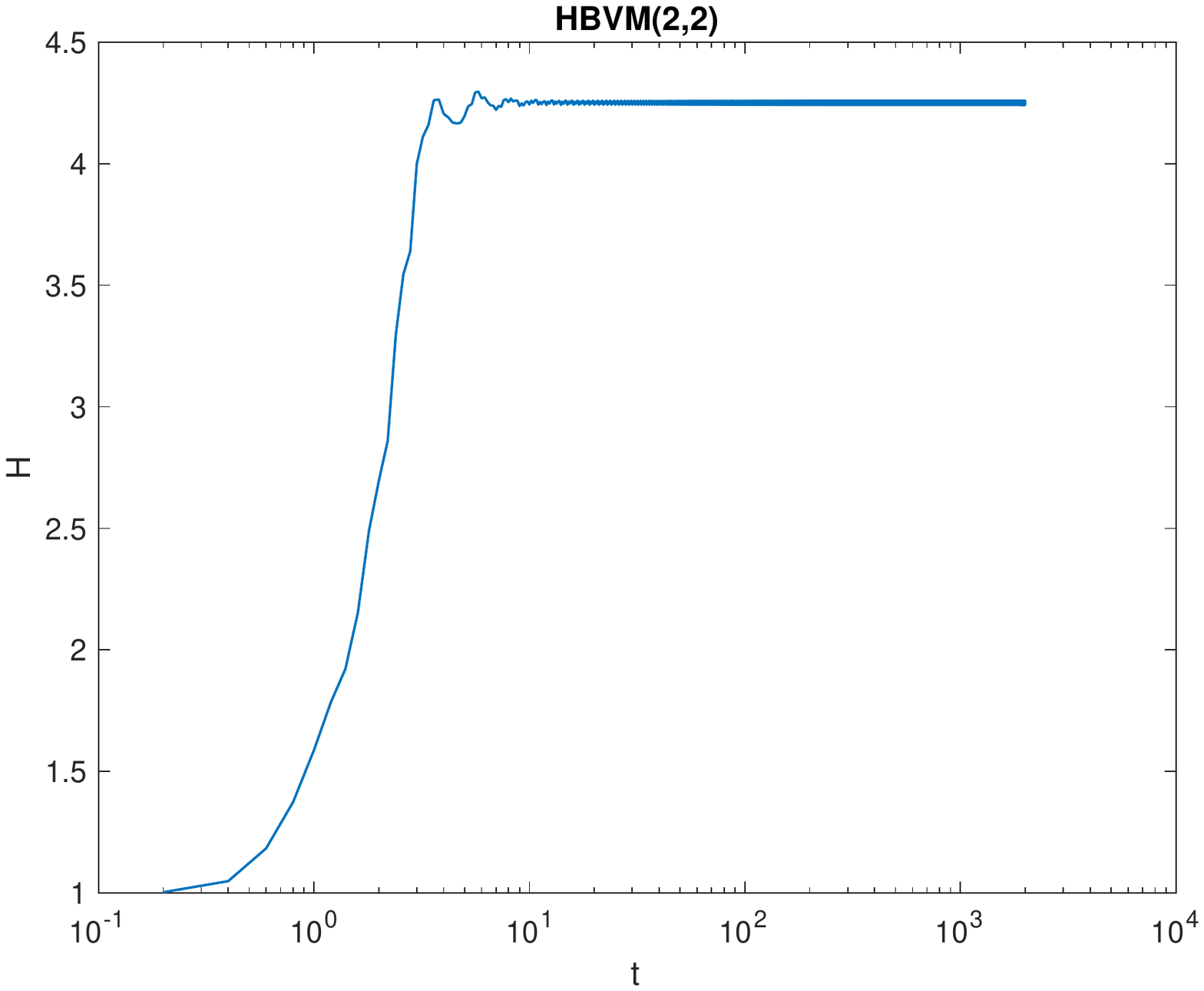}  \includegraphics[width=.55\textwidth]{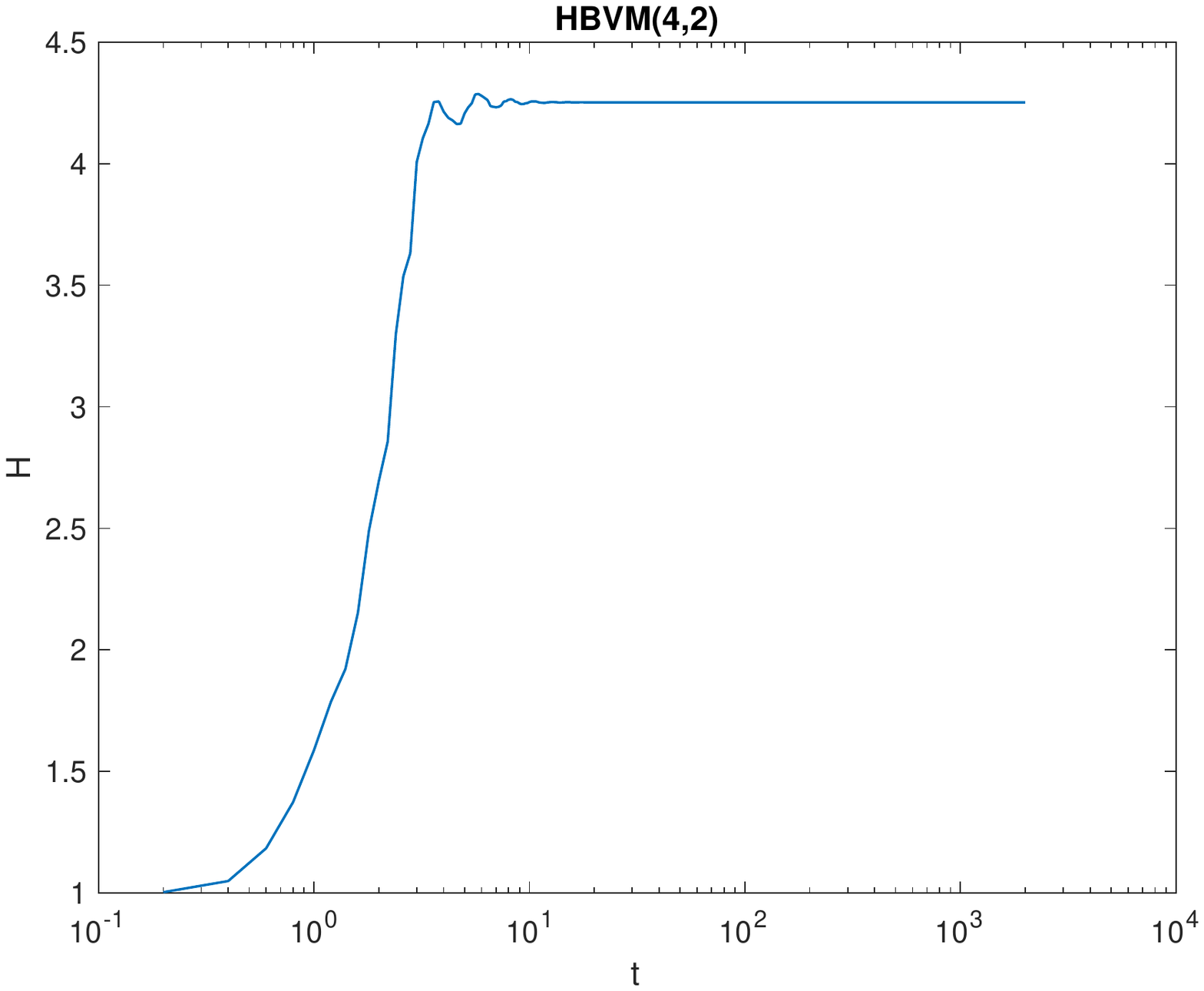}}
\vspace{-5.8cm}

\centerline{ \includegraphics[width=.55\textwidth]{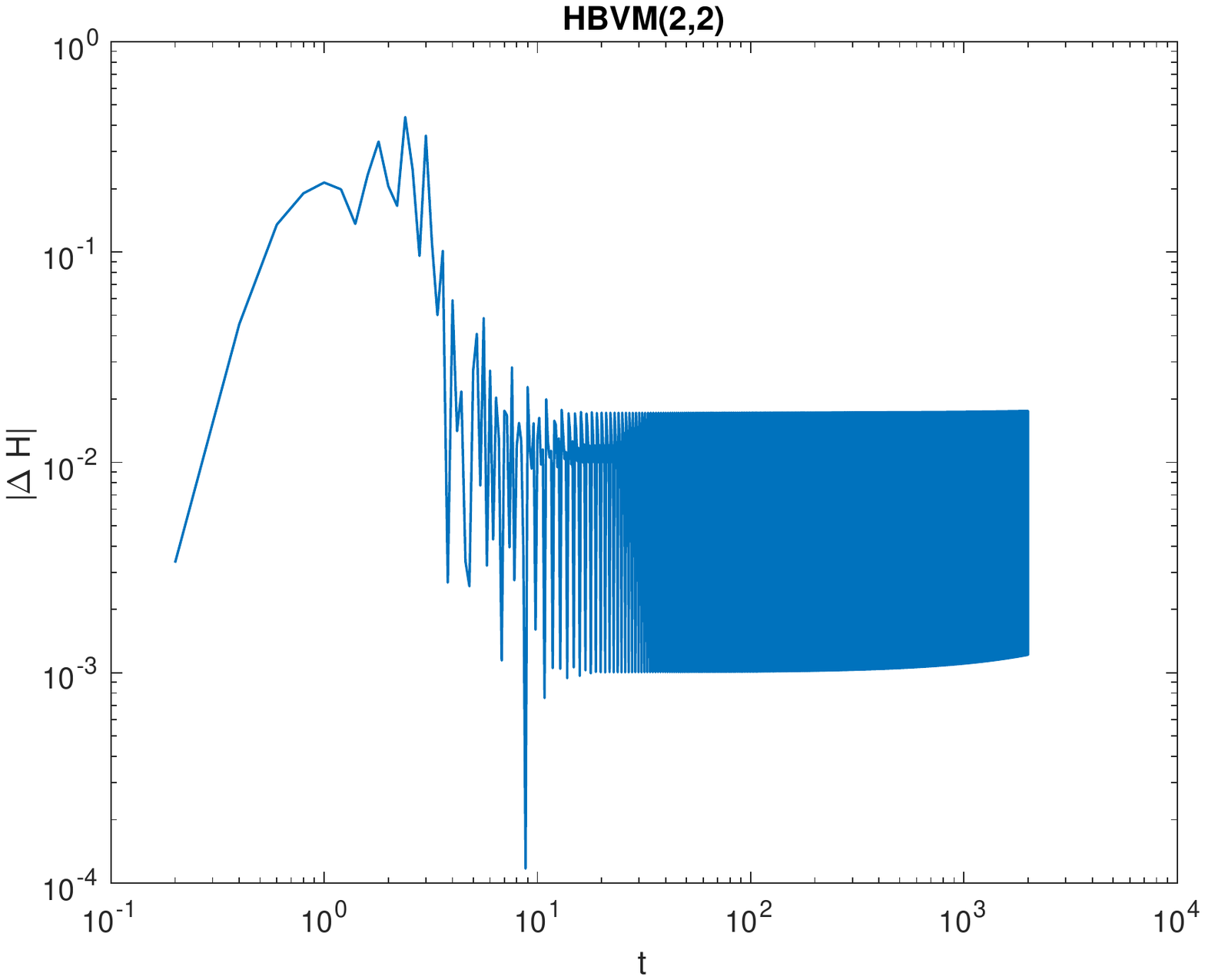}  \includegraphics[width=.55\textwidth]{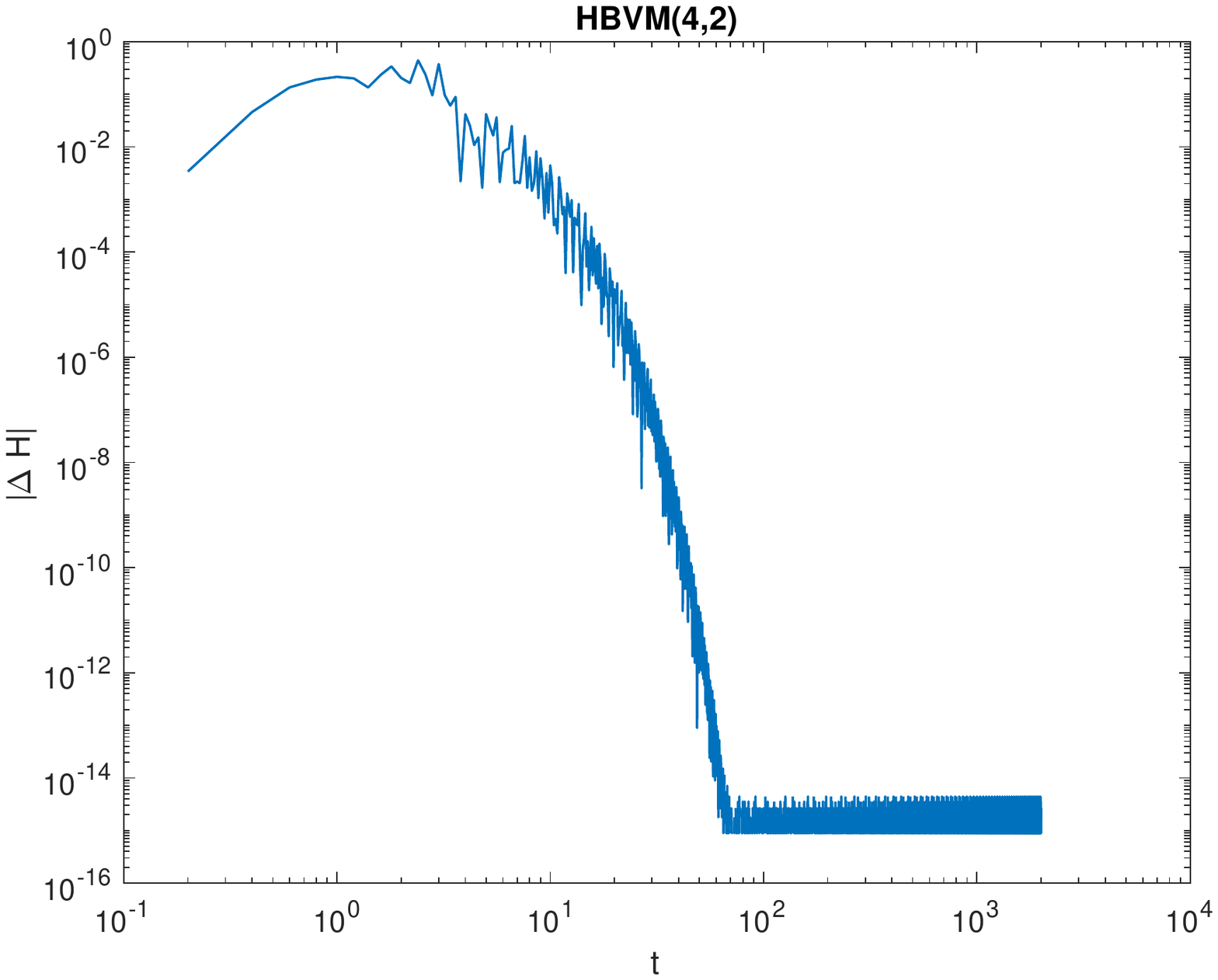}}
\vspace{-5.8cm}

\centerline{ \includegraphics[width=.55\textwidth]{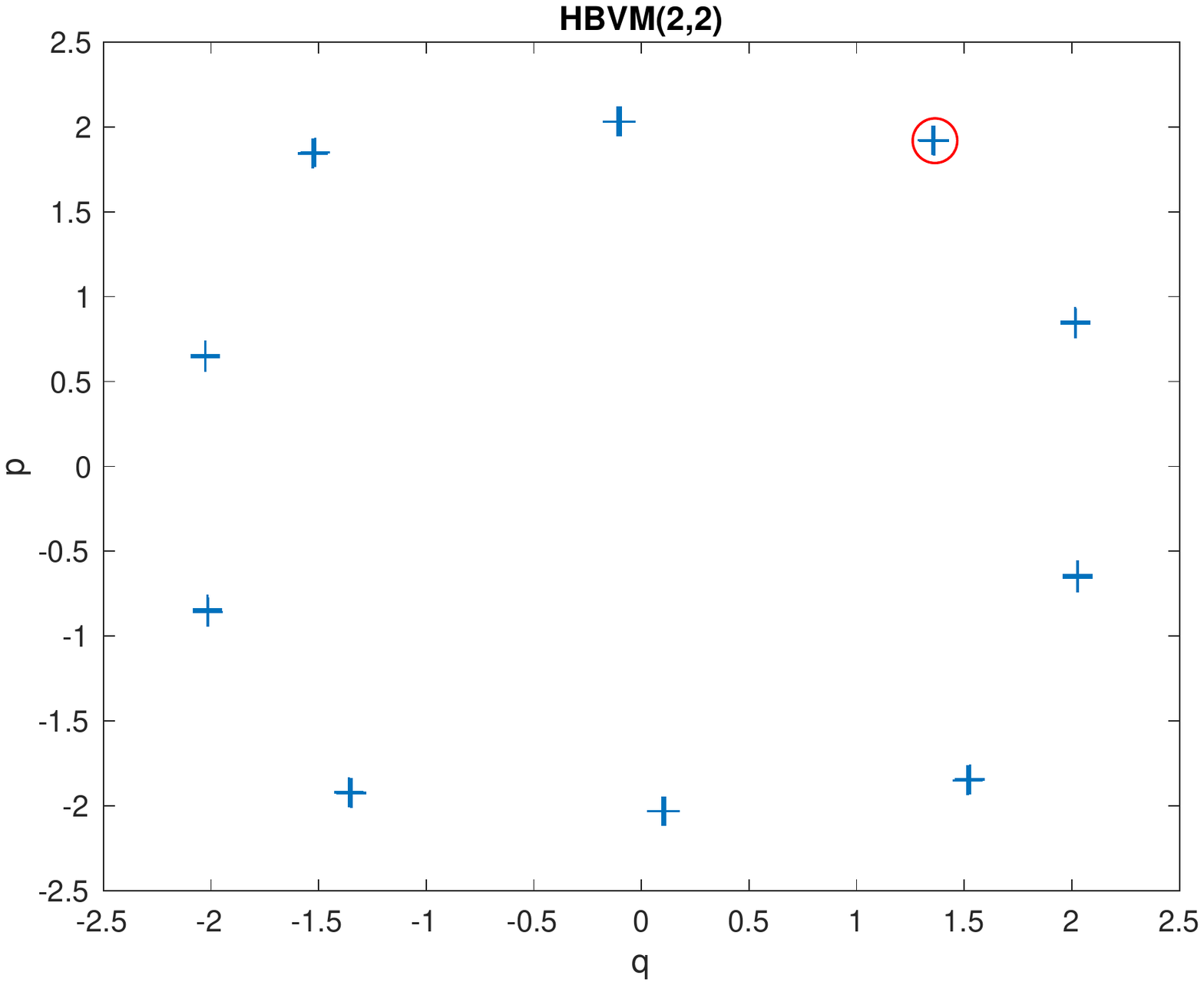}  \includegraphics[width=.55\textwidth]{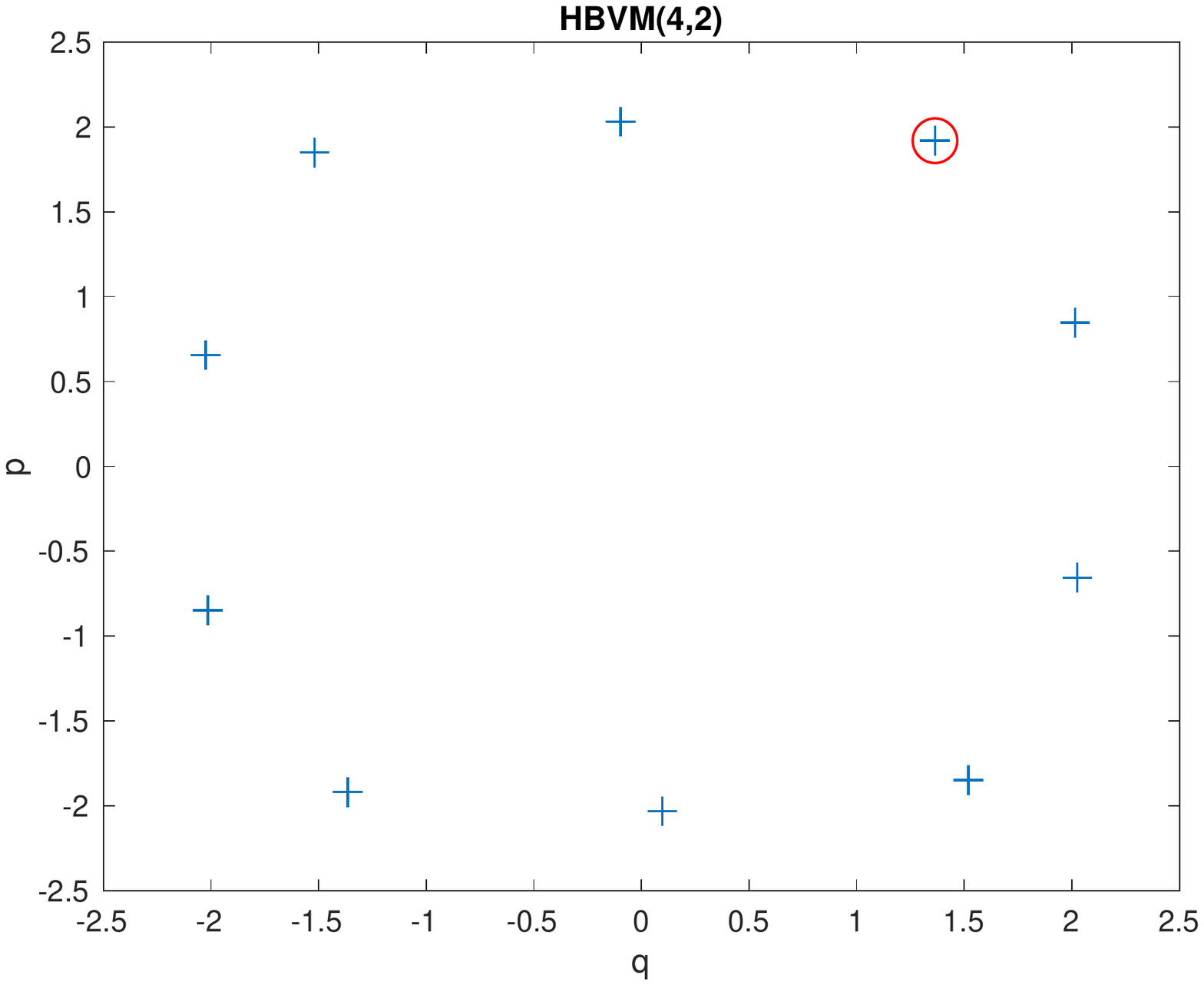}} 
\vspace{-2.8cm}

  \caption{Numerical results for problem (\ref{eq:prob1}) solved by using HBVM(2,2), left plots, and HBVM(4,2), right plots, using a timestep $h=0.2$ (see the text for details).}
  \label{fig:fig1}
\end{figure}

With reference to (\ref{eq:Hdde})--(\ref{eq:inidde}), the first problem is defined as follows:
\begin{equation}\label{eq:prob1}
m=1, \qquad H(q,p) = \frac{1}4\left( q^4 + p^4 \right), \qquad \alpha = 10^{-1},\qquad \tau=1, \qquad \phi(t)\equiv \sqrt 2, \qquad \psi(t) \equiv 0.
\end{equation}
We solve this problem by using the following methods:
\begin{itemize}
\item HBVM(2,2) (i.e., the 2-stage Gauss method),
\item HBVM(4,2).
\end{itemize}
Both methods are fourth-order, according to Theorem~\ref{thm:oden_1}, with HBVM(4,2) energy-conserving once a periodic orbit of integer period is eventually reached (see Theorem~\ref{thm:codipro_1}). Problem (\ref{eq:prob1}) possesses an attracting periodic orbit with period $T=2\tau=2$ which suggests using a stepsize $h$ equal to a submultiple of $\tau$, in order to mimic a corresponding discrete periodic solution.  As we are going to see, unlike the 2-stage Gauss collocation method, the conservation property of HBVM(4,2) results in a precise resolution of this task.
We solve the problem on the interval $[0,2\cdot 10^3]$ by using a timestep $h=\tau/5 = 0.2$.
Figure~\ref{fig:fig1} summarizes the obtained results.
\begin{itemize}
\item In the upper row of the figure are the plots of the numerical Hamiltonian, $H(q_n,p_n)$, from which one deduces that both methods quite soon reach a stationary behavior. 

\item  To better discern the asymptotic behavior of the two numerical solutions, the central pictures show the plots of $\left|H(q_n,p_n)-H(q_{n-1},p_{n-1})\right|$ for the two methods. From these plots one infers that, while the stationary value of the Hamiltonian is constant for the HBVM(4,2) method, it is only approximately constant for the HBVM(2,2) method, with oscillations having amplitude of order $10^{-2}$.

\item The bottom row contains the plots of the numerical trajectory in the phase plane for both methods, relative to the interval $[2\cdot 10^2,2\cdot 10^3]$ (i.e., after the transient phase). For both methods, the solution seems to repeat every 10 points (i.e., with period $T=2\tau=2$). However, the points obtained by the HBVM(2,2) method are not actually periodic, whereas they are  (within to machine precision and independently of the used stepsize $h$) for the HBVM(4,2) method. To confirm this, in Table~\ref{tab:tab1} we list the last 20 points of the trajectories  computed after each period $T=10h$ and lying inside the two small circles highlighted in the plots. As one may see, only the first 4 digit of the points of the trajectory computed by the HBVM(2,2) method are retained, whereas the points computed by the HBVM(4,2) method differ at most on the last digit.

\end{itemize}

\begin{table}[t]
\centering
 \caption{The last 20 points of the trajectories inside the circles in the plots on the bottom line of Figure~\ref{fig:fig1}.}
  \label{tab:tab1}
  
\small
\smallskip
\begin{tabular}{|rr|rr|}
\hline
\multicolumn{2}{|c|}{HBVM(2,2)} & \multicolumn{2}{c|}{HBVM(4,2)}\\
\hline
\multicolumn{1}{|c}{$q$} &\multicolumn{1}{c}{$p$} &\multicolumn{1}{|c}{$q$} &\multicolumn{1}{c|}{$p$}\\
\hline
     1.344913051657652~&     1.924341608176171~&     1.364023296679203~&     1.918490612087558~\\
     1.344895222079097~&     1.924347768593204~&     1.364023296679201~&     1.918490612087558~\\
     1.344877390115245~&     1.924353929533504~&     1.364023296679201~&     1.918490612087558~\\
     1.344859555766308~&     1.924360090996891~&     1.364023296679201~&     1.918490612087558~\\
     1.344841719032502~&     1.924366252983184~&     1.364023296679200~&     1.918490612087558~\\
     1.344823879914032~&     1.924372415492206~&     1.364023296679201~&     1.918490612087558~\\
     1.344806038411116~&     1.924378578523775~&     1.364023296679200~&     1.918490612087558~\\
     1.344788194523964~&     1.924384742077714~&     1.364023296679201~&     1.918490612087559~\\
     1.344770348252789~&     1.924390906153841~&     1.364023296679202~&     1.918490612087558~\\
     1.344752499597800~&     1.924397070751980~&     1.364023296679202~&     1.918490612087558~\\
     1.344734648559217~&     1.924403235871946~&     1.364023296679200~&     1.918490612087559~\\
     1.344716795137255~&     1.924409401513561~&     1.364023296679199~&     1.918490612087559~\\
     1.344698939332121~&     1.924415567676645~&     1.364023296679200~&     1.918490612087559~\\
     1.344681081144034~&     1.924421734361018~&     1.364023296679199~&     1.918490612087559~\\
     1.344663220573212~&     1.924427901566499~&     1.364023296679200~&     1.918490612087558~\\
     1.344645357619866~&     1.924434069292906~&     1.364023296679200~&     1.918490612087558~\\
     1.344627492284208~&     1.924440237540062~&     1.364023296679200~&     1.918490612087558~\\
     1.344609624566458~&     1.924446406307785~&     1.364023296679200~&     1.918490612087559~\\
     1.344591754466830~&     1.924452575595894~&     1.364023296679201~&     1.918490612087558~\\
     1.344573881985545~&     1.924458745404208~&     1.364023296679203~&     1.918490612087558~\\
     \hline
     \end{tabular}
\end{table}     

\subsubsection*{Problem~2}  % dhbvm.m + Htaum1.m
The second example is similar in nature to the previous one but considers a non-polynomial Hamiltonian function with two degrees of freedom. With reference to (\ref{eq:Hdde})--(\ref{eq:inidde}), it is defined by: 
\begin{eqnarray}\label{eq:prob2}
&&m=2, \qquad H(q,p) = \frac{1}4\left( q_1^4 +q_2^4+ p_1^4 +p_2^4\right) + \frac{\pi}2\left(\frac{1}{\|q\|_2^2} + \frac{2}{\|p\|_2^2}  \right),  \\[1mm]
&&\alpha = 5\cdot 10^{-2},\qquad \tau=1, \qquad \phi(t)\equiv (0.1,~1)^\top, \qquad \psi(t) \equiv (1,~0.2)^\top.\nonumber
\end{eqnarray}
Again, we have experienced the existence of a periodic orbit with period $T=2\tau=2$. We solve this problem on the interval $[0,10^3]$ with timestep $h=\tau/10 = 0.1$, by using the following methods:
\begin{itemize}
\item HBVM(2,2) (i.e., the 2-stage Gauss method),
\item HBVM(10,2).
\end{itemize}
Both methods are fourth-order, the latter being {\em practically} energy-conserving, for the given timestep, in the event that a periodic orbit is reached.

\begin{figure}[hp]
\centering

\vspace{-2.5cm}
\centerline{ \includegraphics[width=.55\textwidth]{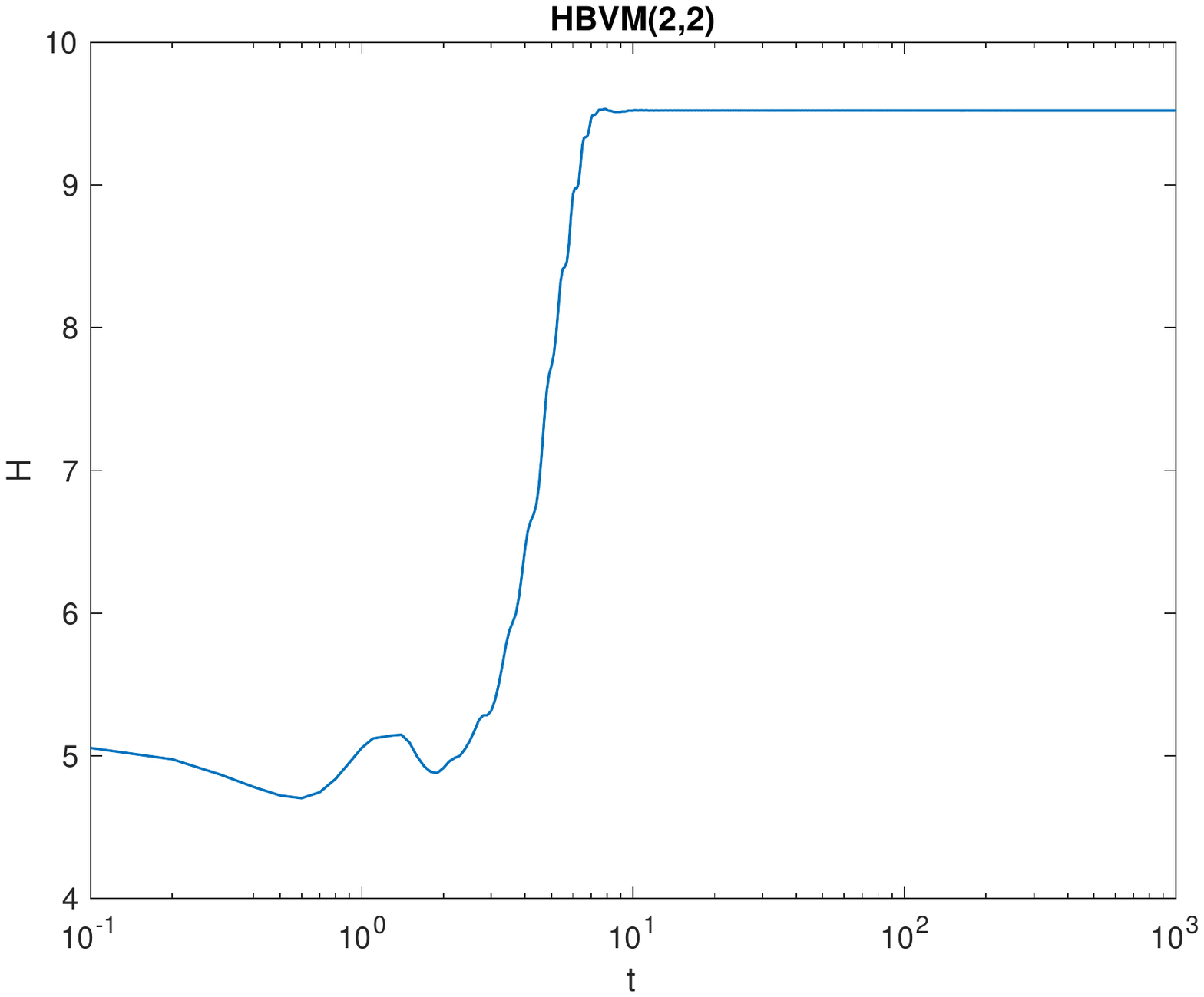}  \includegraphics[width=.55\textwidth]{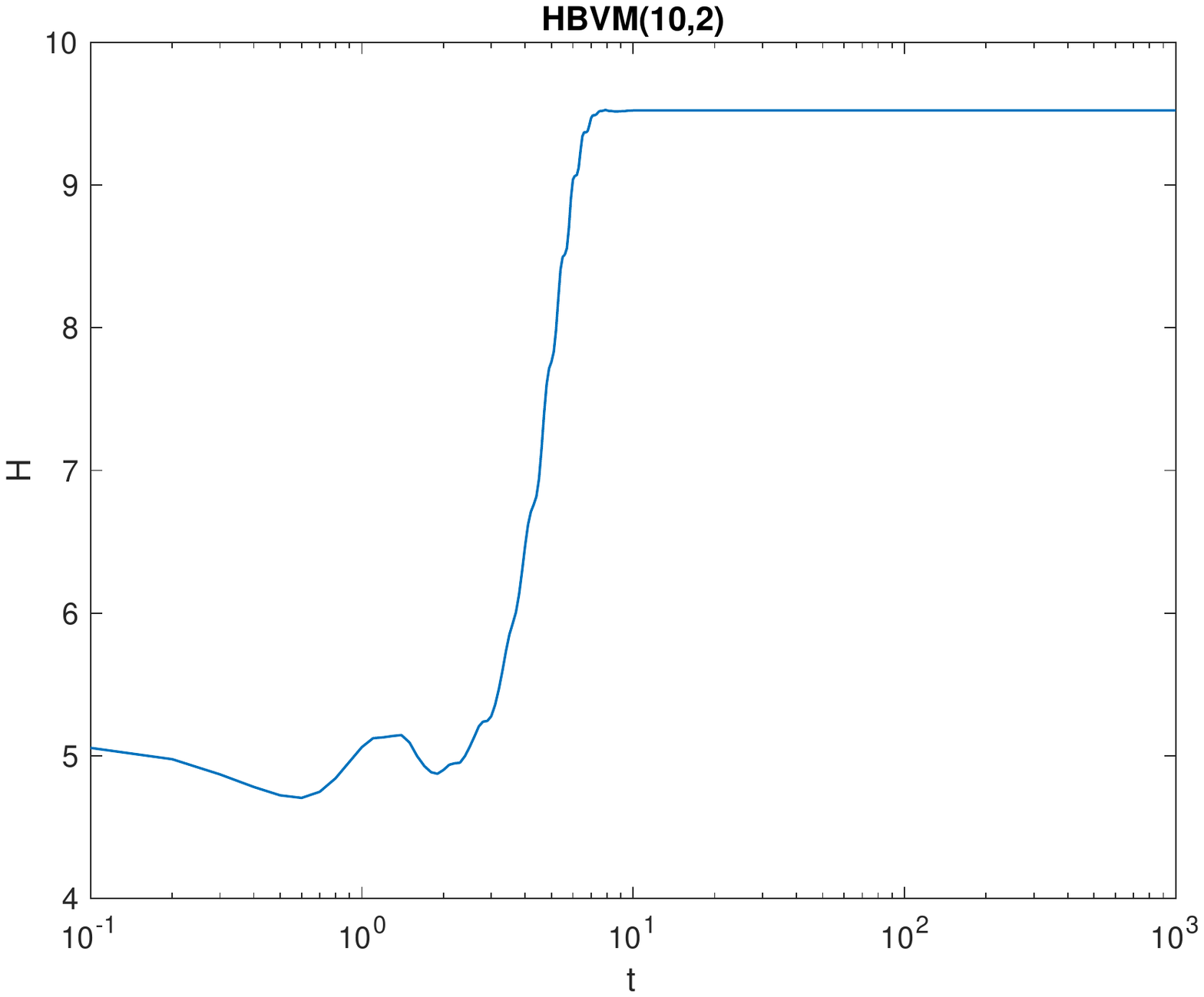}}
\vspace{-5.8cm}

\centerline{ \includegraphics[width=.55\textwidth]{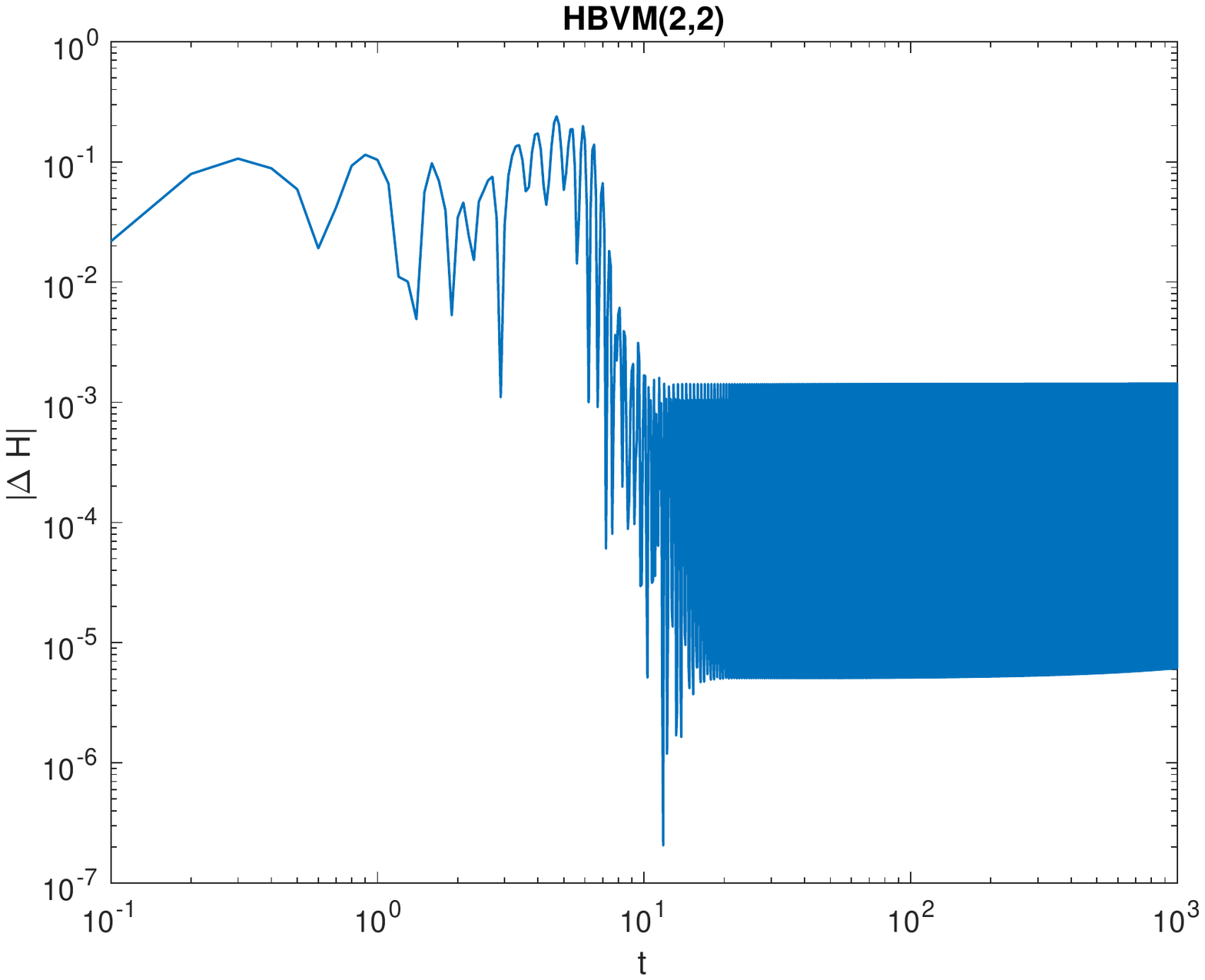}  \includegraphics[width=.55\textwidth]{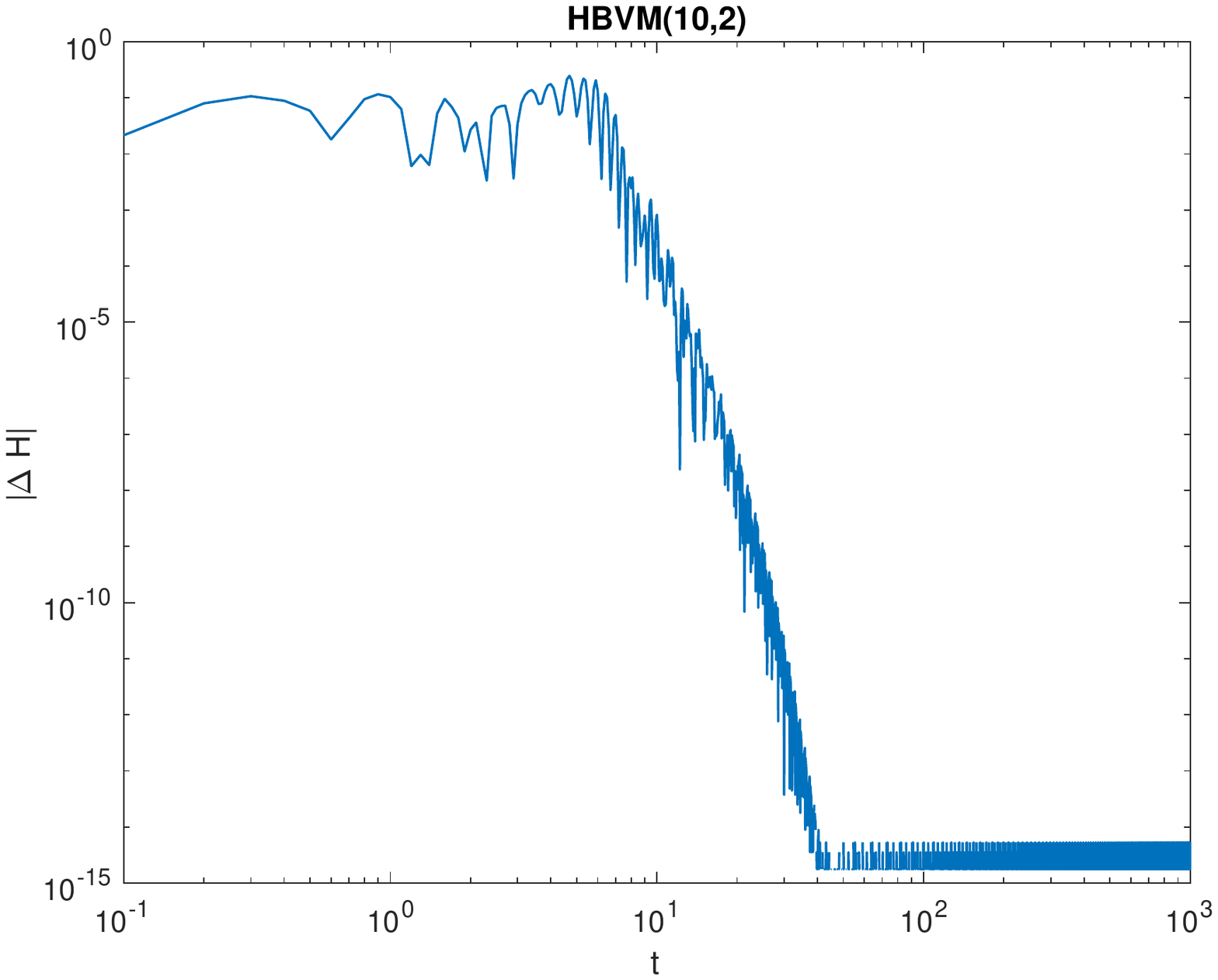}}
\vspace{-5.8cm}

\centerline{ \includegraphics[width=.55\textwidth]{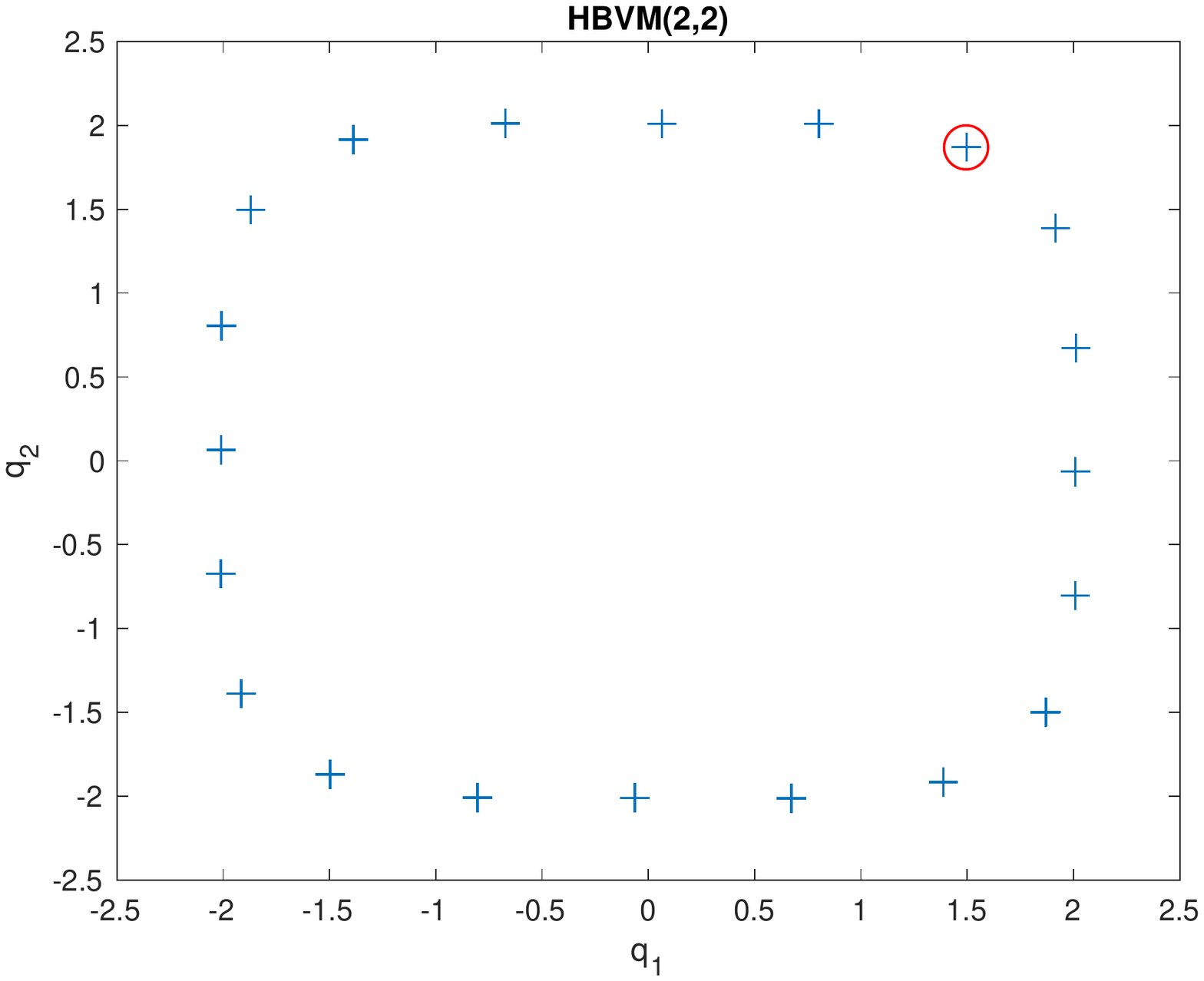}  \includegraphics[width=.55\textwidth]{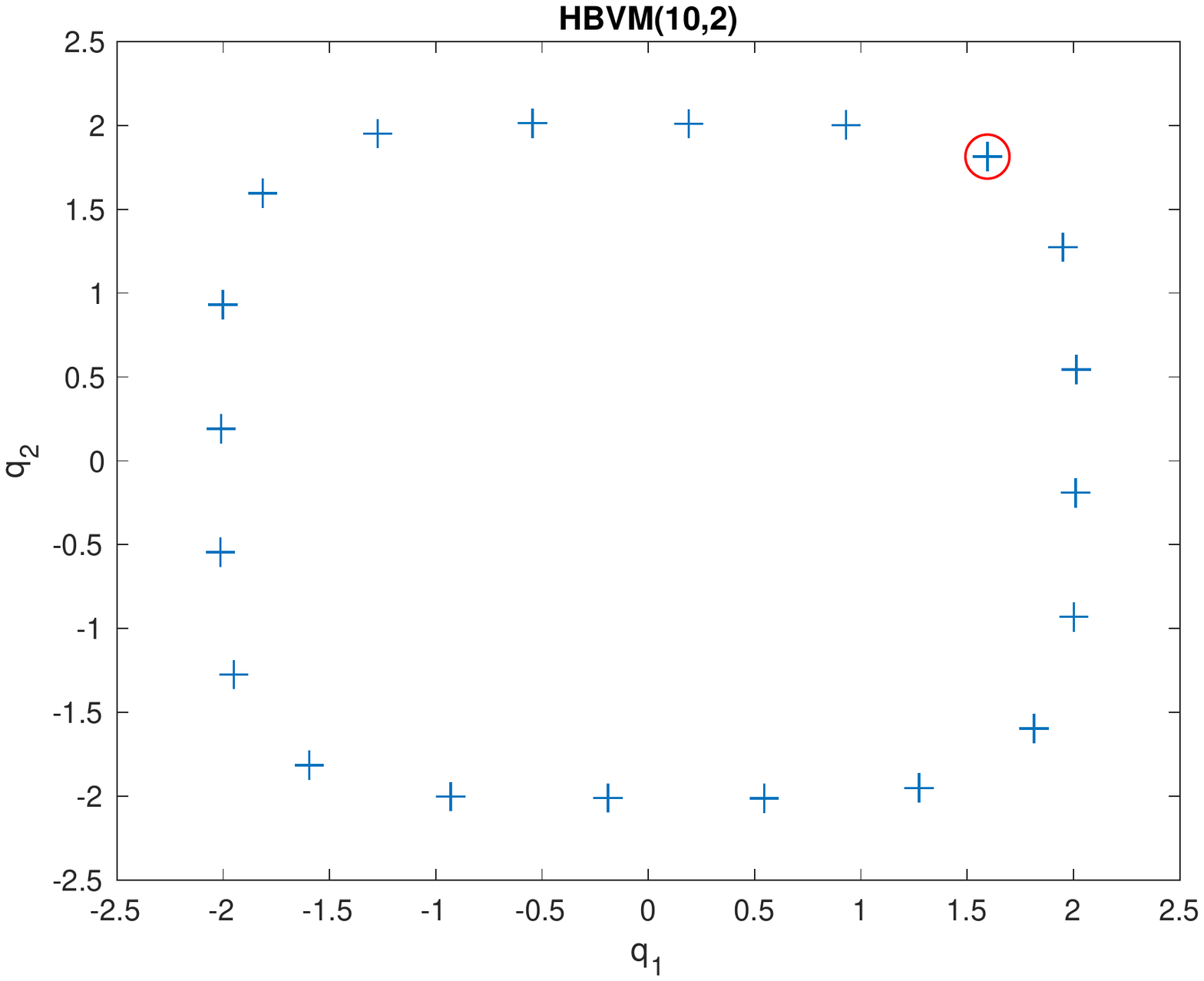}} 
\vspace{-2.8cm}

  \caption{Numerical results for problem (\ref{eq:prob2}) solved by using HBVM(2,2), left plots, and HBVM(10,2), right plots, using a timestep $h=0.1$ (see the text for details).}
  \label{fig:fig2}
\end{figure}

Also in this case, the conservation property of HBVM(10,2)  turns out to be crucial in reproducing a discrete orbit with period precisely equal to 2, while a small phase drift affects the solution yielded by the 2-stage Gauss collocation method. Figure~\ref{fig:fig2}, which is similar to Figure~\ref{fig:fig1}, summarizes the obtained results.
\begin{itemize}
\item In the upper row of the figure are the plots of the numerical Hamiltonian, namely $H(q_n,p_n)$: for both methods it seems to  reach a stationary behavior. 

\item The second row shows the plots of $\left|H(q_n,p_n)-H(q_{n-1},p_{n-1})\right|$ for the two methods. From these plots one infers that, while the stationary value of the Hamiltonian is constant (up to round-off) for the HBVM(10,2) method, it is only approximately constant for the HBVM(2,2) method, with oscillations having amplitude of order $10^{-3}$. 

\item The bottom row contains the plots of the numerical trajectory in the $q_1-q_2$ plane for both methods, relative to the interval $[10^2,10^3]$ (i.e., after the transient phase). For both methods, the solution seems to repeat every 20 points (i.e., with period $T=2\tau=2$). However,  only the points obtained by the HBVM(10,2) method are actually periodic. To confirm this, in Table~\ref{tab:tab2} we list the last 20 points of the trajectories  computed after each period $T=20h$ and lying inside the two small circles displayed in the plots. As one may see, the Gauss collocation method only retain the first 5 digits after each period, whereas the points computed by the HBVM(10,2) method differ at most on the last digit.

\end{itemize}

\begin{table}[t]
\centering
 \caption{The last 20 points of the trajectories inside the circles in the plots on the bottom line of Figure~\ref{fig:fig2}.}
  \label{tab:tab2}
  
\small
\smallskip
\begin{tabular}{|rr|rr|}
\hline
\multicolumn{2}{|c|}{HBVM(2,2)} & \multicolumn{2}{c|}{HBVM(10,2)}\\
\hline
\multicolumn{1}{|c}{$q_1$} &\multicolumn{1}{c}{$q_2$} &\multicolumn{1}{|c}{$q_1$} &\multicolumn{1}{c|}{$q_2$}\\
\hline
   1.500006047618583 &  1.868403720200248  & 1.595245320422993 & 1.813631211153069\\
   1.500014079966090  & 1.868399733165645  & 1.595245320422992  & 1.813631211153067\\
   1.500022113284264  & 1.868395745553114  & 1.595245320422991  & 1.813631211153069\\
   1.500030147573164  & 1.868391757362593  & 1.595245320422994 &  1.813631211153068\\
   1.500038182832831  & 1.868387768594017  & 1.595245320422991  & 1.813631211153067\\
   1.500046219063319  & 1.868383779247324  & 1.595245320422993  & 1.813631211153069\\
   1.500054256264677  & 1.868379789322453  & 1.595245320422994  & 1.813631211153069\\
   1.500062294436967  & 1.868375798819335  & 1.595245320422991  & 1.813631211153067\\
   1.500070333580236  & 1.868371807737912  & 1.595245320422991  & 1.813631211153068\\
   1.500078373694533  & 1.868367816078119  & 1.595245320422993  & 1.813631211153069\\
   1.500086414779919  & 1.868363823839889  & 1.595245320422992  & 1.813631211153067\\
   1.500094456836429  & 1.868359831023165  & 1.595245320422991  & 1.813631211153069\\
   1.500102499864130  & 1.868355837627883  & 1.595245320422994  & 1.813631211153068\\
   1.500110543863066  & 1.868351843653977  & 1.595245320422991  & 1.813631211153067\\
   1.500118588833287  & 1.868347849101385  & 1.595245320422993  & 1.813631211153069\\
   1.500126634774853  & 1.868343853970044  & 1.595245320422994  & 1.813631211153069\\
   1.500134681687813  & 1.868339858259893  & 1.595245320422991  & 1.813631211153067\\
   1.500142729572213  & 1.868335861970863  & 1.595245320422991  & 1.813631211153068\\
   1.500150778428111  & 1.868331865102895  & 1.595245320422993  & 1.813631211153069\\
   1.500158828255558  & 1.868327867655927  & 1.595245320422992  & 1.813631211153067\\
   1.500166879054607  & 1.868323869629889  & 1.595245320422991  & 1.813631211153069\\
   1.500174930825299  & 1.868319871024725  & 1.595245320422994  & 1.813631211153068  \\
      \hline
     \end{tabular}
\end{table}     

\subsubsection*{Problem~3}  % dhbvm.m + Hpend.m

For the last problem, we are no more interested in periodic trajectories. Instead, we consider a {\em delay Hamiltonian problem} with {\em dissipation}. This can be achieved by choosing a negative value of the parameter $\alpha$ in (\ref{eq:qpdde}). With reference to (\ref{eq:Hdde})--(\ref{eq:inidde}), the selected parameters are:
\begin{equation}\label{eq:prob3}
m=1, \quad H(q,p) = \frac{1}2 p^2-\cos q, 
\quad \alpha = -10^{-5},\quad \tau=1, \quad \phi(t)\equiv 0, \quad \psi(t) \equiv 1.99999.
\end{equation}
This problem is a dissipative delay-variant of the nonlinear pendulum, with the initial condition chosen close to the separatrix  (the level set $H(q,p)=1$)  between the two different regimes of the pendulum: librations around the straight-down stationary position, and rotations. For the given initial conditions, the pendulum should undergo damped oscillations  with a decreasing trend of the Hamiltonian function $H(q_n,p_n)$. Consequently, when using relatively large stepsizes, it is fundamental to reproduce the correct dissipation of the Hamiltonian along the numerical trajectory.

\begin{figure}[hp]
\centering

\vspace{-2.5cm}
\centerline{ \includegraphics[width=.55\textwidth]{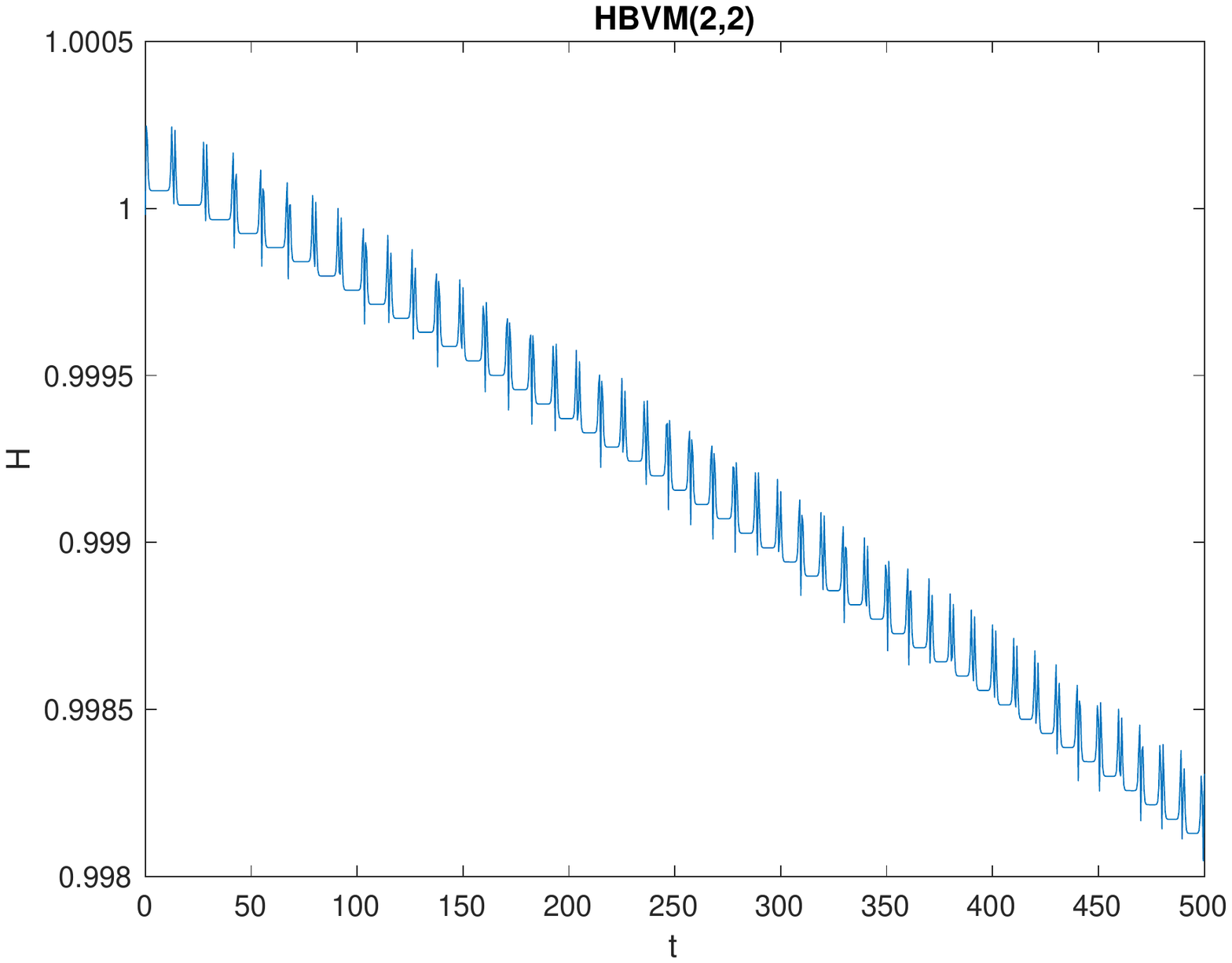}  \includegraphics[width=.55\textwidth]{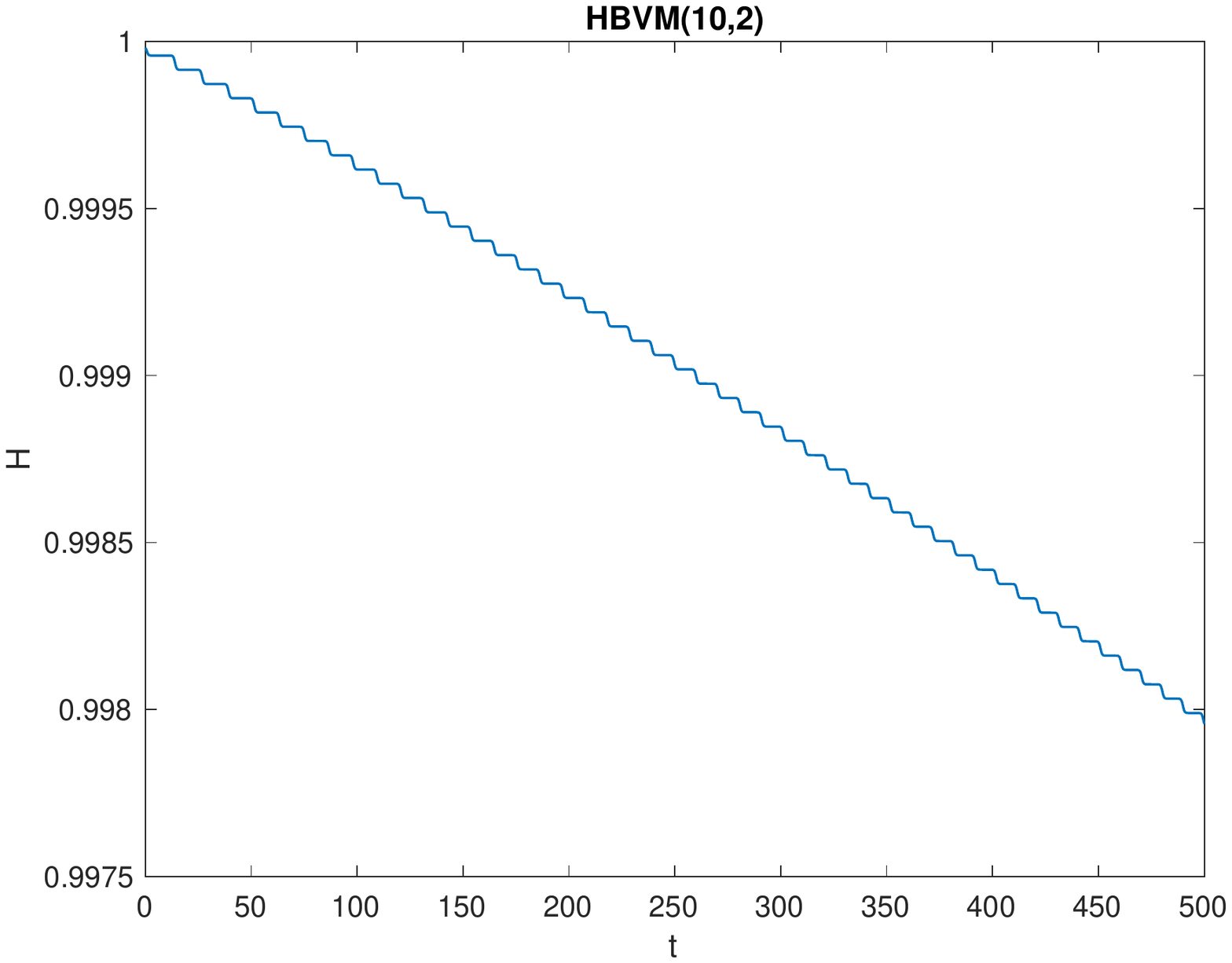}}
\vspace{-5.8cm}

\centerline{ \includegraphics[width=.55\textwidth]{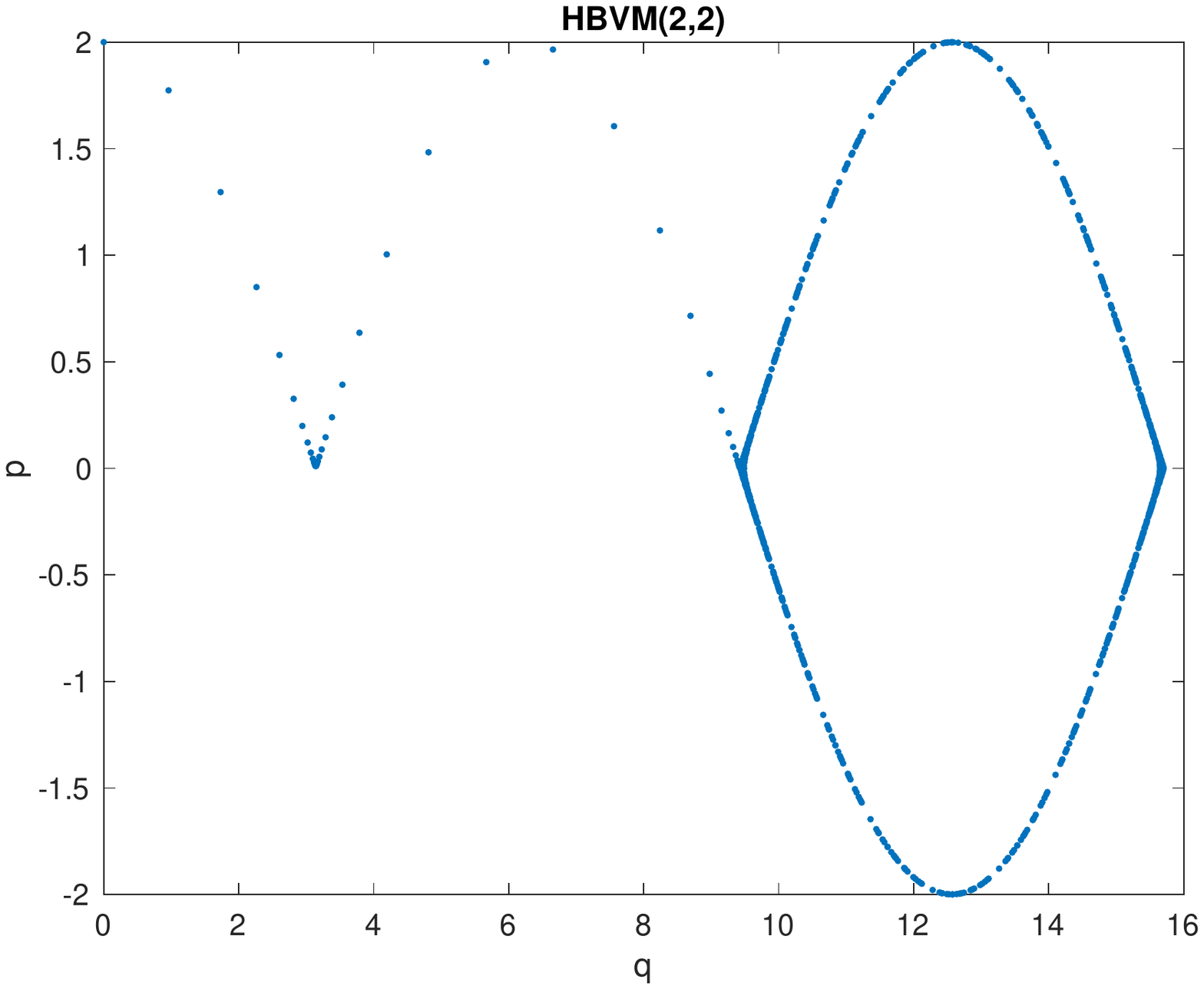}  \includegraphics[width=.55\textwidth]{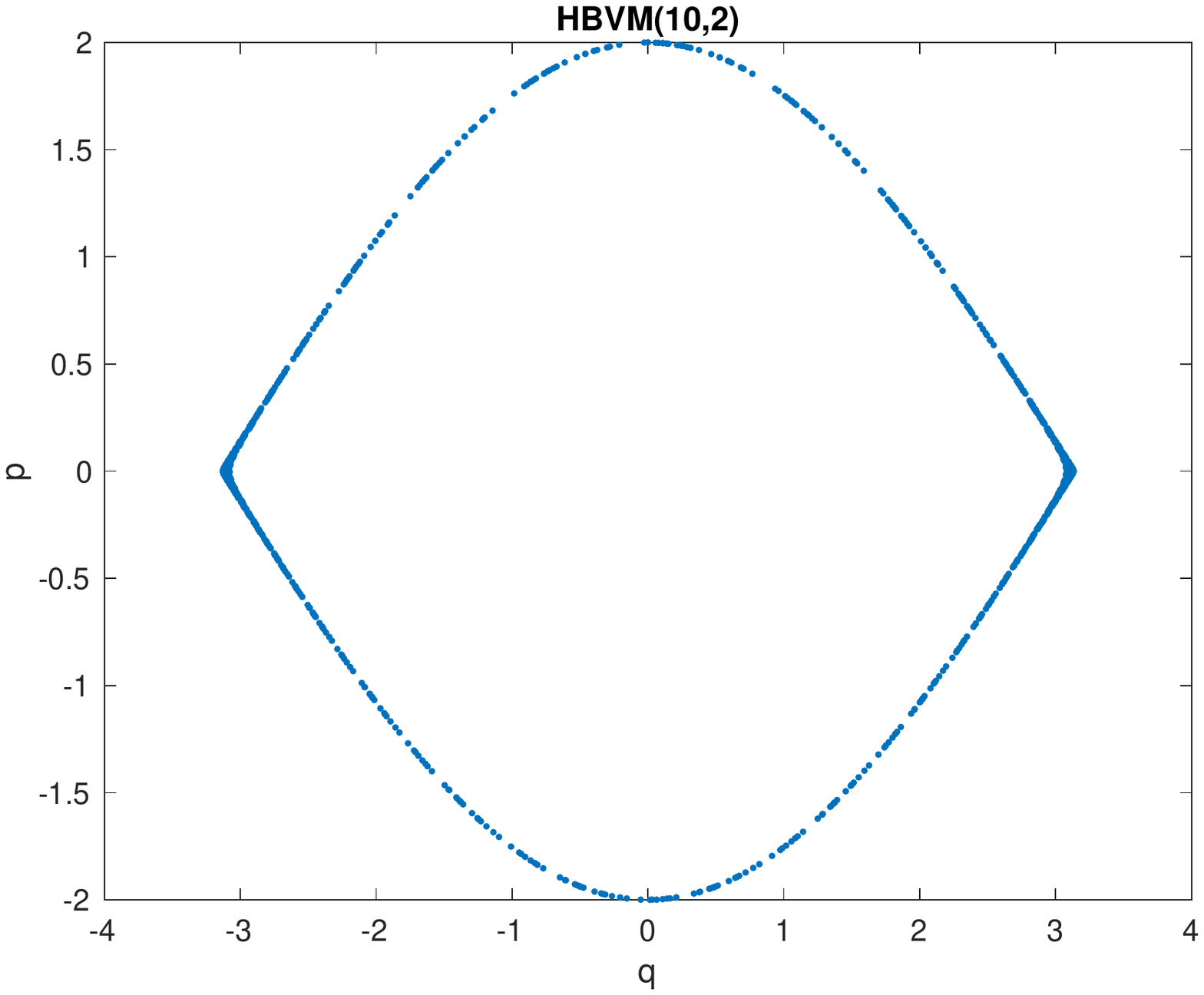}}
\vspace{-5.8cm}

\centerline{ \includegraphics[width=.55\textwidth]{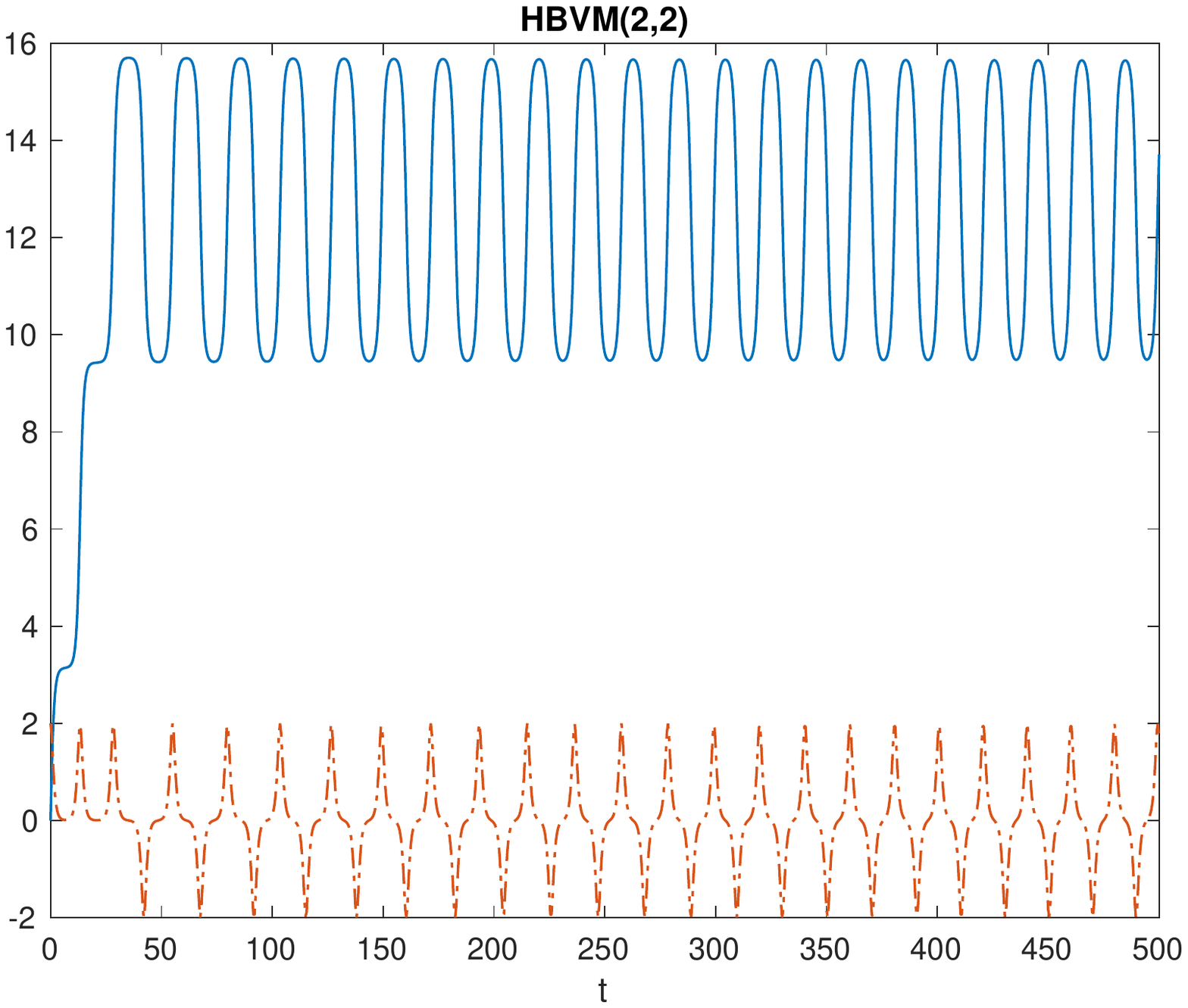}  \includegraphics[width=.55\textwidth]{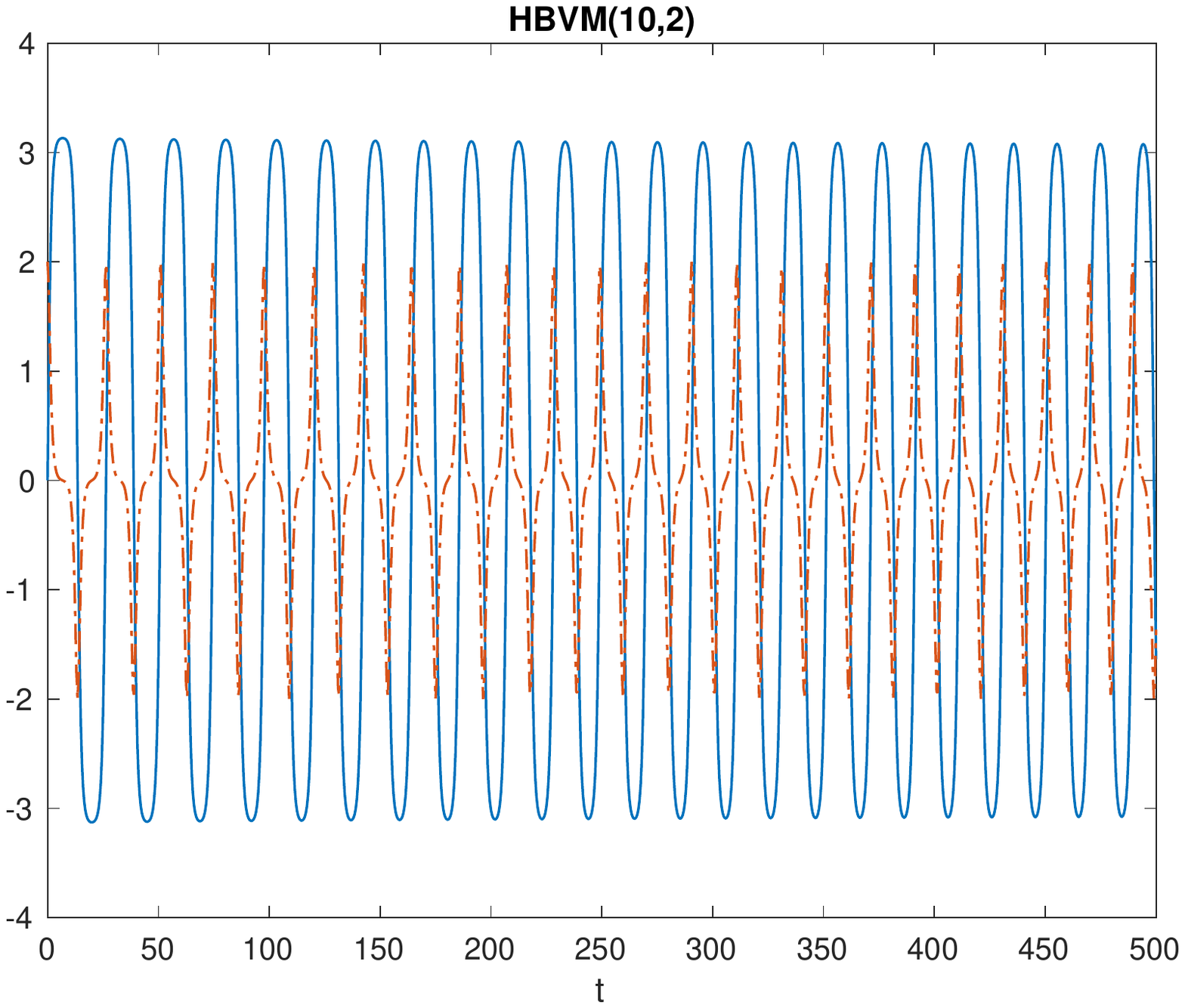}} 
\vspace{-2.8cm}

  \caption{Numerical results for problem (\ref{eq:prob3}) solved by using HBVM(2,2), left plots, and HBVM(10,2), right plots, using a timestep $h=0.5$ (see the text for details).}
  \label{fig:fig3}
\end{figure}

We solve this problem on the interval $[0,500]$, with a timestep $h=\tau/2 = 0.5$, by using the following methods:
\begin{itemize}
\item HBVM(2,2) (i.e., the 2-stage Gauss method),
\item HBVM(10,2).
\end{itemize}
Figure~\ref{fig:fig3} summarizes the obtained results.
\begin{itemize}
\item In the upper row of the figure are the plots of the numerical Hamiltonian, $H(q_n,p_n)$, from which one deduces that both methods have a dissipation trend of the energy $H$. Nevertheless, for HBVM(2,2) the values of the Hamiltonian becomes quite larger than 1 in the initial part of the trajectory  and undergoes fictitious oscillations which cause  the numerical solution  to escape the correct region of the phase space  where the dynamics should take place, as we are going to see.  This is not the case for the  HBVM(10,2) method, whose numerical Hamiltonian  decreases  in the correct way, thus remaining always  smaller than 1.

\item The central pictures show the numerical solution in the phase space. As one may see, the numerical solution provided by HBVM(2,2) ``jumps'' twice, before being trapped into an invariant region.  This means that the pendulum undergoes two complete rotations until it looses enough energy and begins  oscillating around the rest position.  On the contrary, the numerical solution obtained by using HBVM(10,2) always remains in the correct region.

\item The bottom row contains the plots of the numerical solution w.r.t. time, confirming that the numerical solution provided by the HBVM(2,2) method ``jumps'' twice, whereas that obtained by the HBVM(10,2) method does not. 

\end{itemize}

\section{Conclusions}\label{sec:conclusions}
In this paper we have fully developed a thorough approach for obtaining polynomial approximations to the solution of initial value ODE and DDE problems.
It allows us to derive a wide class of Runge-Kutta methods, whose properties are easily discussed within the framework, as well as their actual implementation. Some numerical tests, concerning the numerical simulation of solutions of certain DDE problems of Hamiltonian type, confirm this.
 The present approach  leaves room for generalizations along several directions: in particular to different kind of problems, besides the ones considered here. Another relevant direction of investigation consists in looking for approximations belonging to functional subspaces different than polynomials: that is, by considering orthonormal functional bases different from (\ref{eq:leg}).  Both directions will be the subject of future investigations.

\section*{Funding and conflicts of interests}
The authors have no affiliation with any organization with a direct or indirect financial interest in the subject matter discussed in the manuscript.  The authors acknowledge the financial support from the {\em mrSIR} crowdfunding \cite{mrsir}.

%\section*{Acknowledgments}

\bibliographystyle{siamplain}

\begin{thebibliography}{99}

\bibitem{ABI2019} {\sc P.\,Amodio, L.\,Brugnano, and F.\,Iavernaro.} {\em A note on the continuous-stage Runge-Kutta-(Nystr\"om) formulation of Hamiltonian Boundary Value Methods (HBVMs).} Appl. Math. Comput., 363 (2019) 124634. \url{https://doi.org/10.1016/j.amc.2019.124634} 

\bibitem{ABI2022} {\sc P.\,Amodio, L.\,Brugnano, and F.\,Iavernaro.} {\em Continuous-Stage Runge-Kutta Approximation to Differential Problems.} Axioms, 11 (2022)  192. \url{https://doi.org/10.3390/axioms11050192}

\bibitem{ABI2020} {\sc P.\,Amodio, L.\,Brugnano, and F.\,Iavernaro.} {\em Analysis of Spectral Hamiltonian Boundary Value Methods (SHBVMs) for the numerical solution of ODE problems.} Numer. Algorithms, 83 (2020) 1489--1508. \url{https://doi.org/10.1007/s11075-019-00733-7}

\bibitem{B84} {\sc A.\,Bellen.} {\em One step collocation for delay differential equations.} J. Comput. Appl. Math., 10 (1984) 275--283. \url{https://doi.org/10.1016/0377-0427(84)90039-6}

\bibitem{BZ03}  {\sc A.\,Bellen and M.\,Zennaro.} {\em Numerical Methods for Delay Differential Equations.} Clarendon Press, Oxford, 2003.

\bibitem{BS00}  {\sc P.\,Betsch and P.\,Steinmann.} {\em Conservation properties of a time FE method. I. Time-stepping schemes for $N$-body problems.} Internat. J. Numer. Methods Engrg., 49 (2000) 599--638. \url{https://doi.org/10.1002/1097-0207(20001020)49:5<599::AID-NME960>3.0.CO;2-9}

\bibitem{Bo97}  {\sc C.L.\,Bottasso.} {\em A new look at finite elements in time: a variational interpretation of Runge-Kutta methods.} Appl. Numer. Math., 25 (1997) 355--368. \url{https://doi.org/10.1016/S0168-9274(97)00072-X}

\bibitem{BFCI2014} {\sc  L.\,Brugnano, G.\,Frasca-Caccia, and F.\,Iavernaro.} {\em Efficient  implementation of Gauss collocation and Hamiltonian Boundary Value Methods.} Numer. Algorithms, 65 (2014) 633--650. \url{http://doi.org/10.1007/s11075-014-9825-0} 

\bibitem{BI2016} {\sc L.\,Brugnano and F.\,Iavernaro.} {\em Line Integral Methods for Conservative Problems.}  Chapman et Hall/CRC, Boca Raton, FL, 2016. \url{https://doi.org/10.1201/b19319}

\bibitem{BI2018} {\sc L.\,Brugnano and F.\,Iavernaro.}  {\em Line Integral Solution of Differential Problems.} Axioms, 7(2) (2018) 36. \url{https://doi.org/10.3390/axioms7020036}

\bibitem{BIMR2019} {\sc L.\,Brugnano, F.\,Iavernaro, J.I.\,Montijano, and L.\,R\'andez.} {\em Spectrally accurate space-time solution of Hamiltonian PDEs.} Numer. Algorithms,  81 (2019) 1183--1202. \url{https://doi.org/10.1007/s11075-018-0586-z} 

\bibitem{BIT2010}  {\sc L.\,Brugnano, F.\,Iavernaro, and D.\,Trigiante.} {\em Hamiltonian Boundary Value Methods (Energy Preserving Discrete Line Integral Methods).} JNAIAM J. Numer. Anal. Ind. Appl. Math., 5, no.\,1-2 (2010) 17--37.

\bibitem{BIT2011}  {\sc  L.\,Brugnano, F.\,Iavernaro, and D.\,Trigiante.}  {\em A note on the efficient implementation of Hamiltonian BVMs.} J. Comput. Appl. Math.,  236 (2011) 375--383. \url{https://doi.org/10.1016/j.cam.2011.07.022}

\bibitem{BIT2012} {\sc L.\,Brugnano, F.\,Iavernaro, and D.\,Trigiante.}  {\em A simple framework for the derivation and analysis of effective one-step methods for ODEs.} Appl. Math. Comput., 218 (2012) 8475--8485. \url{https://doi.org/10.1016/j.amc.2012.01.074}

\bibitem{BIT2015} {\sc  L.\,Brugnano, F.\,Iavernaro, and D.\,Trigiante.}  {\em Analisys of Hamiltonian Boundary Value Methods (HBVMs): a class of energy-preserving Runge-Kutta methods for the numerical solution of polynomial Hamiltonian systems.} Commun. Nonlinear Sci. Numer. Simul.,  20 (2015) 650-667. \url{https://doi.org/10.1016/j.cnsns.2014.05.030}

\bibitem{BIZ2021} {\sc  L.\,Brugnano, F.\,Iavernaro, and P.\,Zanzottera.}  {\em A multiregional extension of the SIR model, with application to the COVID-19 spread in Italy.} Math. Meth. Appl. Sci., 44 (2021) 4414--4427. \url{https://doi.org/10.1002/mma.7039}


\bibitem{BM2002} {\sc L.\,Brugnano and C. Magherini.} {\em Blended implementation of block implicit methods for ODEs.} Appl. Numer. Math., 42 (2002) 29--45. \url{https://doi.org/10.1016/S0168-9274(01)00140-4}

\bibitem{BM2009} {\sc L.\,Brugnano and C. Magherini.} {\em Recent advances in linear analysis of convergence for splittings for solving ODE problems.} Appl. Numer. Math.,  59 (2009) 542--557. \url{https://doi.org/10.1016/j.apnum.2008.03.008}

\bibitem{BMR2019} {\sc L.\,Brugnano, J.I.\,Montijano, and L.\,R\'andez.} {\em On the effectiveness of spectral methods for the numerical solution of multi-frequency highly-oscillatory Hamiltonian problems.} Numer. Algorithms, 81 (2019) 345--376. \url{http://dx.doi.org/10.1007/s11075-018-0552-9}

\bibitem{Br04} {\sc H.\,Brunner}. {\em Collocation Methods for Volterra Integral and Related Functional Equations.} Cambridge University Press, Cambridge, 2004.

\bibitem{CMcLMcLOQW2009} {\sc E.\,Celledoni, R.I.\,McLachlan, D.\,McLaren, B.\,Owren, G.R.W.\,Quispel, and W.M.\,Wright.} {\em Energy preserving Runge-Kutta methods.}  M2AN Math. Model. Numer. Anal., 43 (2009) 645--649. \url{https://doi.org/10.1051/m2an/2009020}

\bibitem{DB2008} {\sc G.\,Dahlquist and \AA.\,Bj\"ork.} {\em Numerical Methods in Scientific Computing, Volume I.} SIAM, Philadelphia, 2008.

\bibitem{DRB1986} {\sc J.G.\, Dos Reis and R.L. Baroni.} {On the existence of periodic solutions for autonomous retarded functional-differential equations on $R^2$.}  Proc. Roy. Soc. Edinburgh Sect., A 102 (1986) 259--262. \url{https://doi.org/10.1017/S0308210500026342}

\bibitem{ED2011} {\sc E.\,Engel and R.M.\,Dreizler.} {\em Density Functional Theory, an advanced course.} Springer, Berlin, 2011.

\bibitem{FM2010} {\sc D.\,Furihata and T. Matsuo.} {\em Discrete Variational Derivative Method: A Structure-Preserving Numerical Method for Partial Differential Equations.} Chapman and Hall/CRC, Boca Raton, FL, 2010.

\bibitem{Ha2010}  {\sc E.\,Hairer.} {\em Energy-preserving variants of collocation methods.} JNAIAM J. Numer. Anal. Ind. Appl. Math., 5, no.\,1-2 (2010) 73--84.

\bibitem{HNW08} {\sc E.\,Hairer, S.P.\,N\o rsett, and G.\,Wanner.} {\em Solving Ordinary Differential Equations I, nonstiff problems. Second revised edition (3rd printing).} Springer, Heidelberg, 2008.

\bibitem{HW2002} {\sc E.\,Hairer and G.\,Wanner.} {\em Solving Ordinary Differential Equations I, nonstiff problems. Second revised edition.} Springer, Heidelberg, 2002.

\bibitem{Hu72_1}  {\sc B.L.\,Hulme.} {\em One-step piecewise polynomial Galerkin methods for initial value problems.} Math. Comp., 26 (1972) 415--426. \url{https://doi.org/10.1090/S0025-5718-1972-0321301-2}

\bibitem{Hu72_2}  {\sc B.L.\,Hulme.} {\em Discrete Galerkin and related one-step methods for ordinary differential equations.} Math. Comp., 26 (1972) 881--891. \url{https://doi.org/10.1090/S0025-5718-1972-0315899-8}

\bibitem{IP2007} {\sc F.\,Iavernaro and  B.\,Pace.} {\em $s$-Stage trapezoidal methods for the conservation of Hamiltonian functions of polynomial type.} AIP Conf. Proc. 936 (2007) 603--606. \url{https://doi.org/10.1063/1.2790219}

\bibitem{IP2008} {\sc F.\,Iavernaro and B.\,Pace.} {\em Conservative block-boundary value methods for the solution of polynomial Hamiltonian systems.} AIP Conf. Proc. 1048 (2008) 888--891. \url{https://doi.org/10.1063/1.2991075}

\bibitem{IT2009} {\sc F.\,Iavernaro, and D.\,Trigiante.} {\em High-order symmetric schemes for the energy conservation of polynomial Hamiltonian problems.} JNAIAM J. Numer. Anal. Ind. Appl. Math., 4, no.\,1-2 (2009) 87--111.

\bibitem{KY1974} {\sc J.L.\,Kaplan and J.A.\,Yorke.} {\em Ordinary differential equations which yield periodic solutions of differential delay equations.} J. Math. Anal. Appl., 48 (1974) 317--324. \url{https://doi.org/10.1016/0022-247X(74)90162-0}

\bibitem{PaNu2011} {\sc J.\,Mallet-Paret and R.D.\,Nussbaum}. {\em Stability of periodic solutions of state-dependent delay-differential equations}.  J. Differential Equations, 250 (2011) 4085--4103. \url{https://doi.org/10.1016/j.jde.2010.10.023}

\bibitem{MB2016} {\sc Y.\,Miyatake and J.C.\,Butcher.} {\em A characterization of energy-preserving methods and the construction of parallel integrators for Hamiltonian systems.} SIAM J. Numer. Anal., 54 (2016) 1993--2013. \url{https://doi.org/10.1137/15M1020861}

\bibitem{Nu1973} {\sc R.D.\,Nussbaum}. {\em Periodic solutions of some nonlinear, autonomous functional differential equations.}   Bull. Amer. Math. Soc., 79 (1973) 811--814. \url{https://doi.org/10.1016/0022-0396(73)90053-3}

\bibitem{Nu1973_1} {\sc R.D.\,Nussbaum}. {\em Periodic solutions of some nonlinear, autonomous functional differential equations. II}.  J. Differential Equations, 14 (1973) 360--394. \url{https://doi.org/10.1090/S0002-9904-1973-13330-0}
 
 \bibitem{Nu1979} {\sc R.D.\,Nussbaum}. {\em Uniqueness and nonuniqueness for periodic solutions of $x'(t) = -g(x(t-1))$}. 
 J. Differential Equations, 34 (1979) 25--54. \url{https://doi.org/10.1016/0022-0396(79)90016-0}

\bibitem{QMcL2008} {\sc G.R.W.\,Quispel and D.I.\,McLaren.} {\em A new class of energy-preserving numerical integration methods.} J. Phys. A, 41 (2008) 045206. \url{https://doi.org/10.1088/1751-8113/41/4/045206}
 
 \bibitem{W1975} {\sc H.-O.\,Walther}. {\em Existence of a non-constant periodic solution of a nonlinear autonomous functional differential equation representing the growth of a single species population.}  J. Math. Biol., 1 (1975) 227--240. \url{https://doi.org/10.1007/BF01273745}

\bibitem{mrsir} \url{https://www.mrsir.it/en/about-us/}

\end{thebibliography}

\end{document}